\documentclass
{ar2e}

\usepackage{epsfig,latexsym,xspace,oupbib,hyperref,url}

\newcommand{\ie}{{\it i.e.\/}\xspace}
\newcommand{\eg}{{\it e.g.\/}\xspace}

\newcommand{\viz}{{\it viz.\/}\xspace}

\newcommand{\HIDE}[1]{ }
\newcommand{\COMMENT}[1]{ }

\renewcommand{\dag}{\mbox{$\mathcal D$}}

\newcommand{\half}{\frac{1}{2}}

\newcommand{\eqref}[1]{\mbox{(\ref{eq:#1})}}

\newcommand{\secref}[1]{\mbox{\S$\,$\ref{sec:#1}}}

\newcommand{\figref}[1]{\mbox{Figure~\ref{fig:#1}}}

\newcommand{\defref}[1]{\mbox{Definition~\ref{def:#1}}}

\newcommand{\itref}[1]{\mbox{\ref{it:#1}}}

\newcommand{\exref}[1]{\mbox{Example~\ref{ex:#1}}}

\newcommand{\parents}[1]{{\rm pa}(#1)}

\newcommand{\nondesc}[1]{{\rm nd}(#1)}

\newcommand{\nd}{\nondesc}

\newcommand{\cd}{\,|\,}

\newcommand{\cov}{{\rm Cov}}

\newcommand{\E}{{\mbox{\rm E}}}

\newcommand{\cip}{\mbox{$\perp\!\!\!\perp$}}

\newtheorem{expl}{Example}
\newtheorem{rem}{Remark}
\newtheorem{definer}{Definition}
\newtheorem{theorem}{Theorem}
\newtheorem{algor}{Algorithm}
\newtheorem{thm}[theorem]{Theorem}

\newtheorem{regime}{Regime}

\newtheorem{assumption}{Assumption}

\newtheorem{regimes2}{Regime}

\newcommand{\halm}{\hspace*{\fill} $\Box$\par}

\newenvironment{ex}{\begin{expl}\rm}{\halm\end{expl}}
\newenvironment{defn}{\begin{definer}\rm}{\halm\end{definer}}

 \newcommand{\indo}[2]{\mbox{$#1 \,\cip\,
    #2$}} \newcommand{\ind}[3]{\mbox{$#1 \, \cip\, #2 \cd
    #3$}}\newcommand{\inda}[4]{\mbox{$#1 \perp_{\mbox{\scriptsize
        $#4$}} #2 \cd #3$}} \newcommand{\nind}[3]{\mbox{$#1
    \not\!\!\cip\, #2 \cd #3$}} 
\newcommand{\indp}[4]{\mbox{$#1 \,\cip\, #2 \cd #3$ $[#4]$}}
 \newtheorem{cond}{Condition}
\newcommand{\condref}[1]{\mbox{Condition~\ref{cond:#1}}}
  \renewcommand{\theenumi}{(\roman{enumi})}

\newcommand{\bY}{\mbox{\boldmath$Y$}}
\newcommand{\bE}{\mbox{\boldmath$E$}}

\newcommand{\bm}{\mbox{\boldmath$\mu$}}
\newcommand{\bo}{\mbox{\boldmath$0$}}
\newcommand{\bSigma}{\mbox{\boldmath$\Sigma$}}
\newcommand{\norm}{\mbox{${\cal N}$}}
\newcommand{\idle}{\mbox{$\emptyset$}}

 \newcommand{\var}{\mbox{\rm
    var}}

\newcommand{\ice}{\mbox{\rm ICE}\xspace}
\newcommand{\ace}{\mbox{\rm ACE}\xspace}
\newcommand{\sce}{\mbox{\rm SCE}\xspace}
\newcommand{\ett}{\mbox{\rm ETT}\xspace}

\begin{document}



\jname{Annual Review of Statistics and Its Applications}
\jyear{2014}
\jvol{2}
\ARinfo{1056-8700/97/0610-00}

\title{Statistical Causality from a Decision-Theoretic Perspective}

\markboth{A.~P.~Dawid}{Statistical Causality}

\author{A.~Philip~Dawid
\affiliation{University of Cambridge, UK}}


\begin{abstract}
  We present an overview of the decision-theoretic framework of
  statistical causality, which is well-suited for formulating and
  solving problems of determining the effects of applied causes.  The
  approach is described in detail, and is related to and contrasted
  with other current formulations, such as structural equation models
  and potential responses.  Topics and applications covered include
  confounding, the effect of treatment on the treated, instrumental
  variables, and dynamic treatment strategies.
\end{abstract}

\maketitle
 
\section{Introduction}
\label{sec:intro}
After decades of neglect, recent years have seen a flowering of the
field of statistical causality, with an impressive array of
developments both theoretical and applied.  However there is as yet no
one fully accepted foundational basis for this enterprise.  Rather,
there is a variety of formal and informal frameworks for framing and
understanding causal questions, and hot discussion about their
relationships, merits and demerits: we might mention, among others,
structural equation modelling and path analysis \cite{wright},
potential response models \cite{dbr:as}, functional models
\cite{pearl:book}, and various forms of graphical representation.
This plethora of ``foundations'' leaves statistical causality in much
the same state of confusion as probability theory before Kolmogorov.

The aim of this overview paper is to add to the confusion by
describing one particular approach, based on decision-theoretic
principles, that its author considers superior to the others (in this
of course he is fully aligned with others' attitudes to their own
works).  I will present the main features of the approach, relate it
to and compare it with some other approaches, and show how it works in
some simple applications.  I consider this approach more
straightforward philosophically and mathematically (it requires only a
very small extension to standard statistical methods), and easier to
comprehend and manipulate, than other approaches that introduce new
ingredients and structures, such as potential responses or
deterministic functional relationships.  While I do not expect
wholesale conversion to ths point of view, I hope that readers already
knowledgeable in causal inference will, at the very least, find it
helpful to look at familiar topics through a fresh pair of spectacles.

For further details and developments of the material in this paper,
see
\textcite{apd:cinfer,apd:infdiags,apd:hsss,vd/apd/sg:uai,apd:kent,apd/vd:uai08,apd:beware,apd:pearl,dd:ss,hg/apd:aistats2010,apd:aberdeen,sgg/apd:ett,apd/vd:princstrat,carlo/phil/luisa:overview,apd:causDTchapter,carlo/phil/vanessa:dynamicchapter,apd/pc:2014}.
The lecture notes \textcite{apd:gtp} contain a fuller exposition of
the decision-theoretic approach.  

\subsection*{Causality and agency}

``Causality'' has been a focus of interest for philosophers for
millennia, but --- as befits any worthwhile philosophical conundrum
--- all this attention has not resulted in a settled approach to
understanding it.  In \textcite{apd:pearl} I reviewed a variety of
philosophical conceptions and theories, focusing particularly on that
which is most germane to my own approach: the {\em agency theory\/} of
causality \cite{price:bjps91,hausman:book,woodward:book,woodward:sep}.
This interprets causality as being all about how an external
manipulation would affect a system: for example, how the quality of
the chemical product would respond to adjustments of the lever that
controls the pressure in the production process.  Much of Statistical
Science---in particular, the whole subfield of Experimental
Design---aims to address exactly these kinds of questions about the
effects of interventions on a system, which are indeed a major object
of all scientific enquiry.

An important advantage of the agency approach is its clear separation
of cause and effect variables, resulting in the elimination of
definitional ambiguities associated with the possibility of reverse
causation or common causes.  But having clean definitions is not
enough for practical purposes: we must be able to relate those
definitions to properties of the empirical world.  Whereas this is
relatively unproblematic in cases involving genuine experimentation,
it becomes a major headache when we can only observe a system in its
natural habitat, and are unable to apply to it the interventions that
are essential to understanding its causal properties.  The importance
of making a clear distinction between intervening and merely observing
has been stressed by numerous authors, including
\textcite{dbr:as,meek:bjps,sgs:book,pearl:book}.  Most of the recent
emphasis of the enterprise of ``statistical causality'' focuses on
observational situations, and attempts to identify conditions under
which it is possible to extract causal conclusions from them, and
techniques for doing so.

\subsection*{Effects of causes and causes of effects}
The emphasis in this work will be entirely on the problem of assessing
the future effects of an intervention in a system: that is, on
identifying the ``effects of causes'' (EoC).  An entirely different
problem is that of judging what might have been the cause of an
observed outcome: that is to say, identifying the ``causes of
effects'' (CoE).  This is very much an issue in legal cases, which may
seek to assign responsibility or blame.

Most of the effort to date in statistical causality has focused on CoE
problems, an important exception being \textcite{pearl:book} who
explores both problems in detail.  However, whereas Pearl, and such
others who have dealt with CoE, have used the identical mathematical
machinery (in Pearl's case, based on assumed functional relationships)
for both EoC and CoE, I do not consider this appropriate.  In
\textcite{apd:cinfer} I discussed the relationships and differences
between these problems, and their requisite infrastructures, in some
detail, and argued that, whereas (unlike for EoC) some form of
counterfactual logic appears unavoidable for assessing CoE, the
difficulties in appropriately modelling a CoE problem have been
underappreciated.  Some discussion and analysis of CoE issues in the
context of using epidemiological data to address a legal case of toxic
tort can be found in
\textcite{apd:aberdeen,sef/dlf/apd:socmeth,apd/mm/sef:ba}.

\subsection*{Plan of paper}

I start in \secref{decprob} with a simple example that locates
statistical causality firmly within the purview of classical
statistical decision analysis.  Section~\ref{sec:altform} then
explores some variant formulations of this problem, including
structural equations and potential responses.  In \secref{expt} these
approaches are explored and compared in the familiar context of
statistical experimental design.  Section~\ref{sec:obsconf} moves the
discussion on to the more problematic context of causal inference from
observational data, and explores the meaning and representation of the
important concept of ``no confounding'': from a decision-theoretic
approach, this is usefully described in terms of relationships between
different ``regimes'' ---\eg, interventional or observational ---
under which data can, in principle at least, be gathered.  It is shown
how this and similar requisite properties can be usefully expressed
and manipulated in terms of an extension of the probabilistic notion
of conditional independence to allow for both stochastic and
non-stochastic variables.  Section~\ref{sec:ci} develops the
associated algebraic and graphical theory.  In \secref{CAUSDAG} we
introduce influence diagrams as useful graphical representations of
causal problems, and relate these to the use of directed acyclic graph
representations as described by \textcite{pearl:book}.  The remainder
of the paper explores, from the decision-theoretic perspective, a
number of important special applications.  Section~\ref{sec:ident}
examines the observational identification of causal effects using
sufficient covariates, propensity analysis and ``{\em
  do\/}-calculus'', and the possibility of identifying ``the effect of
treatment on the treated''; \secref{instrument} considers the use of
instrumental variables; and \secref{dynamic} treats problems where a
sequence of actions can be applied, over time, in response to
intermediate observations and outcomes.  The discussion in
\secref{disc} expresses some scepticisms about currently popular
approaches that make unavoidable use of counterfactual reasoning.

\section{A Decision Problem}
\label{sec:decprob}
I have a headache and am considering whether or not to take two
aspirin tablets.  It is generally accepted that aspirin has a
beneficial effect on headaches: in some sense --- and it will be our
task to try to make this more precise --- taking aspirin {\em
  causes\/} headaches to get better faster.  The solution to my
decision problem is thus intimately bound up with the cause-effect
relationship of aspirin on headache.

This observation leads naturally on to a suspicion that some
part,\footnote{Specifically, the part that aims to understand the
  effects of applied causes (EoC); problems of identifying the causes
  of observed effects (CoE) require a different approach.} at least,
of the enterprise of ``statistical causality'' might be fruitfully
recast as a special application of standard statistical decision
analysis.  However, this point of view is not a currently popular one,
and indeed there is a variety of other attempts to interpret causality
in a statistical setting.  We shall be considering the relationships,
similarities and differences between these, and hope to demonstrate
that the decision-theoretic approach is more natural, more
straightforward, and more useful than its competitors.

To formulate the decision problem, let the binary variable $X$ denote
whether I take aspirin ($X=1$) or not ($X=0$), and let $Y$ be the
log-time it takes for my headache to go away.  I myself can choose
$X$: it is a decision variable, and does not have a probability
distribution.  Nevertheless, it is still meaningful to consider my
conditional distribution, $P_x$ say, for how the eventual response $Y$
will turn out, given that I choose $X=x$.  For the moment we assume
the distributions $P_0$, $P_1$ to be known.  Where we need to be
definite, we shall (purely for simplicity) take them to have the
following normal probability density function:
\begin{equation}
  \label{eq:normal}
  p(y \cd X=x) = (2\pi\sigma^2)^{-\half} \exp -\frac{(y - \mu_x)^2}{2\sigma^2},
\end{equation}
with a mean $\mu_0$ or $\mu_1$ according as $x= 0$ or $1$, and
variance $\sigma^2$ in either case.
  
The distribution $P_1$ [resp., $P_0$] can be interpreted as expressing
my {\em hypothetical\/} uncertainty about $Y$, {\em if\/} I were to
decide on action $X=1$ [resp., $X=0$].  It can incorporate various
sources and types of uncertainty, including stochastic effects of
external influences arising and acting between the points of treatment
application and eventual response.  The distributions $P_1$ and $P_0$
are all that is needed to address my decision problem: I simply need
to compare the two different hypothetical distributions for $Y$,
decide which one I prefer, and take the associated decision.

One possible comparison of $P_1$ and $P_0$ might be in terms their
respective means, $\mu_1$ and $\mu_0$; the ``effect'' of taking
aspirin, rather than nothing, might then be quantified by means of the
change in the expected response, $\delta:= \mu_1 - \mu_0$.
Alternatively, we might look at the difference of the means of $Z =
e^Y$ under the two possible treatments:
$e^{\sigma^2/2}(e^{\mu_1}-e^{\mu_2})$.  Or we might compare the
variances of $Z$ under the two treatments.  Any such comparison of an
appropriately chosen feature of the two hypothetical distributions of
$Y$ can be regarded as a summary of the {\em causal effect\/} of
taking aspirin (as against taking nothing).

More formally, we might apply statistical decision analysis (see \eg\
\textcite{raiffa}) to structure and solve this decision problem.
Suppose that I quantify the loss that I will suffer if my headache
lasts $y$ minutes by means of a real-valued loss function, $L(y)$.  If
I were to take the aspirin, my expected loss would be $E_{Y\sim P_1}
\{L(Y)\}$; if not, it would be $E_{Y\sim P_0} \{L(Y)\}$.  The
principles of statistical decision analysis now direct me to choose
the treatment leading to the smaller expected loss.  A trivial but
fundamentally important point is that, whatever loss function is used,
this solution will only involve the two hypothetical distributions,
$P_1$ and $P_0$, for $Y$, conditional on taking either action.  The
``effect of treatment'' might be measured by the reduction in expected
loss, $E_{P_0} \{L(Y)\} - E_{P_1} \{L(Y)\}$: and the correct decision
will be to take aspirin just when this reduction is positive.

Although there is no uniquely appropriate measure of ``the effect of
treatment'', in the rest of this paper for simplicity we shall focus
on the difference of the means of the two hypothetical distributions:
$\delta:= E_{P_1}(Y) - E_{P_0}(Y)$.

\section{Alternative Formulations}
\label{sec:altform}

\subsection{Decision-theoretic (DT) model}
\label{sec:stochmod}

The essential ingredients of the decision-theoretic analysis above
were the two hypothetical distributions for $Y$, conditional on $X=0$
and on $X=1$.  These were specialised to be normal:
\begin{equation}
  \label{eq:normal2}
  Y \cd X=x  \sim \norm(\mu_x, \sigma^2).
\end{equation}
We term this a {\em stochastic\/} or {\em decision-theoretic\/} (DT)
model.  In this formulation, the term {\em average causal effect\/}
(\ace) simply denotes the difference of the means of the two
hypothetical distributions for $Y$, $\mu_1-\mu_0$.

\subsection{Simple structural equation (SSE) model}
\label{sec:FRAME-sem}
An alternative way in which the assumptions of \eqref{normal2} are
very often expressed is as follows:
\begin{equation}
  \label{eq:struct}
  Y = \mu_X + E
\end{equation}
where
\begin{equation}
  \label{eq:enorm}
  E \sim \norm(0,\sigma^2)
\end{equation}
and it is implicit that the ``error'' $E$ is independent of $X$.  A
system of equations such as \eqref{struct}, which may contain hundreds
of relationships representing response (``endogenous'') variables as
functions of other (both endogenous and ``exogenous'') variables as
well as of external ``error'' variables such as $E$, together with
associated explicit or implicit assumptions about the joint
distribution of the error terms, constitutes a simple {\em Structural
  Equation\/} (SSE) model.  Such models are popular in econometrics
and other fields as representations of causal structures.

The assumptions of the SSE model \eqref{struct} clearly imply the
distributional properties of the DT model \eqref{normal}.  Does this
mean the SSE model is equivalent to the DT model?  To assume this
would be to ignore the additional {\em algebraic\/} structure of
\eqref{struct}, whereby $Y$ is represented as a deterministic
mathematical function of the two variables $X$ and $E$.  Unlike the
distributional formulation of \eqref{normal}, in \eqref{struct} all
the uncertainty is compressed into the single variable $E$, {\em
  via\/} \eqref{enorm}.  If we take \eqref{struct} and its ingredients
seriously, we can get more out of it.

It is common, and indeed seems very natural, implicitly to interpret
\eqref{struct} as follows.  The values of $X$ and $E$ are assigned
separately (by the decision maker and by Nature, respectively), and
$Y$ is then determined by the equation.  In this case, given that $E$
takes value $e$, then if I set $X$ to $x$, $Y$ will take value $y_x:=
\mu_x + e$.  That is, I will observe the variable $Y_x: = \mu_x + E$.
We can regard $Y_x$ as the {\em potential response\/} to the
hypothetical setting $X=x$.  It will become the actual response, $Y$,
when indeed $X=x$: $Y = Y_X = \mu_X + E$.  Note that, in this
interpretation, when in fact I set $X=1$, the {\em
  counterfactual\/}\footnote{since predicated on a hypothesis that
  runs counter to known facts} response\ $Y_0 = \mu_0 + E$ to $X=0$ is
still a well-defined function of the ingredients of model
\eqref{struct}.

With the above interpretations, if I were to switch my decision from
$X=1$ to $X=0$, then $Y$ would switch from $Y_1 = \mu_1 +E$ to $Y_0 =
\mu_0 + E$, with the identical $E$.  The ``causal effect'' of this
switch might then be measured by $Y_1-Y_0$, which in this case is the
constant $\mu_1 - \mu_0$.  This is a purely algebraic comparison,
unrelated to the stochastic properties of the model.  It may be termed
an {\em individual causal effect\/} (\ice).

Note that all these manipulations rely fundamentally on the implicit
assumption that the value of $E$ remains fixed, irrespective of which
decision I take.  Such an assumption simply has no counterpart in the
DT model \eqref{normal2}.

\subsection{Extended structural equation (ESE) model}
\label{sec:extension}
An extension of the SSE model~\eqref{struct} is given by:
\begin{equation}
  \label{eq:estruct}
  Y = \mu_X + E_X
\end{equation}
where we now have a pair $\bE = (E_0, E_1)$ of error variables, having
a bivariate distribution, assumed independent of $X$.  In particular,
when $X=0$ we have $Y = Y_0 := \mu_0 + E_0$, with $E_0$ having its
initially assigned distribution; and similarly $Y = Y_1 := \mu_1 +
E_1$ when $X=1$.  And $Y = Y_X$.

Suppose we model the pair of errors $\bE$ as bivariate normal with
standard normal margins:
\begin{equation}
  \label{eq:smbiv}
  \bE \sim \norm(\bo, \bSigma),
\end{equation}
where $\bSigma$ ($2 \times 2$) has diagonal entries $\sigma^2$ and
off-diagonal entries $\rho\sigma^2$.  Then (no matter what the
correlation $\rho$ may be) the DT model \eqref{normal} for $Y$ given
$X$ will be obtained.  But again, if we take the algebraic structure
seriously, and further suppose that the value of $\bE$ is unaffected
by the choice made for $X$, we can go further and define potential
responses $Y_x := \mu_x + E_x$, as well as the \ice, $Y_1-Y_0$, which
in this case is a random quantity, $\mu_1-\mu_0 + E_1 - E_0$.

It is important to note that the relationship between the ESE model
and the induced DT model is many-one: the dependence structure (here
embodied in the correlation $\rho$) does not enter into the induced DT
model.

\subsection{Potential response (PR) model}
\label{sec:prmod}
As seen above, starting from a structural equation model (whether
simple or extended), we can define the pair $\bY = (Y_0,Y_1)$ of
potential responses, and derive its bivariate distribution, in terms
of the ingredients of that model.  For our ESE model above, the
implied distribution of $\bY$ is
\begin{equation}
  \label{eq:PRbiv}
  \bY \sim \norm(\bm, \bSigma)
\end{equation}
with $\bm := (\mu_0, \mu_1)$.  Both the value and the distribution of
$\bY$ are regarded as independent of the applied treatment $X$.

We might alternatively start at this point, simply taking the pair
$\bY$ as a primitive ingredient of our model, having a bivariate
distribution (for example, that of \eqref{PRbiv}), and again assuming
that both the value and the distribution of $\bY$ will be unchanged if
we change the value of $X$.  This is the general {\em potential
  response\/} (PR) model.  The underlying philosophical conception is
that both potential responses are real and coexist --- even though it
is logically impossible to observe both of them together.

Starting from a PR model we can recover a DT model: in particular, if
we start with \eqref{PRbiv} we recover \eqref{normal2}.  But again,
this relationship is many-one.

\subsection{Functional model}
\label{sec:FRAME-fm}
Mathematically, the models introduced in \secref{FRAME-sem},
\secref{extension} \secref{prmod} all have the following common {\em
  functional\/} form:
\begin{equation}
  \label{eq:func}
  Y = f(X, U),
\end{equation}
where $X$ is a decision variable representing the cause of interest;
$Y$ is the effect of interest; $U$ is a further extraneous random
variable whose value and distribution are taken as independent of $X$;
and $f$ is a deterministic function of its arguments.

In the ESE model~\eqref{estruct}, we can take $U = \bE$ and $f(x,
(e_0, e_1)) = \mu_x + e_x$; the SSE model of \eqref{struct} is the
degenerate case of this having $U = E$ and $f(x, e) = \mu_x + e$.

In the case of a PR model, we can formally take $U$ to be the pair
$(Y_0, Y_1)$, and the function $f$ to be given by:
\begin{equation}
  \label{eq:prf}
  f(x, (y_0, y_1)) = y_x.
\end{equation}

In all the above applications, the variable $U$ typically represents a
somewhat imaginary quantity, that does not correspond to any variable
observable in the empirical world.  Indeed, were that to be made a
requirement, the application of functional models would be limited to
the very special situation of complete determinism, with the pair of
real variables $(X,U)$ fully determining $Y$.

A general functional model of the form \eqref{func} is mathematically
equivalent to a PR model, if we define $Y_0 = f(0, U)$, $Y_1 = f(1,
U)$.  We thus see that (mathematically if not necessarily in terms of
their interpretation) PR models, ESE models and general functional
models need not be distinguished.  Further, any functional model
determines a DT model: under model \eqref{func} the relevant
distribution of $Y$, given $X=x$, is simply the marginal distribution
of $f(x,U)$.  Conversely, given any DT model, $Y \cd X=x \sim P_x$, we
can construct a functional model corresponding to it in this way: one
simple way is as a potential response model \eqref{prf}, in which the
marginal distribution of $Y_x$ is $P_x$.  However, in contrast to the
essentially unique cross-correspondence between the other models
considered above, the functional representation of a DT model is far
from unique---as can again be seen, for example, from the
arbitrariness of the dependence parameter $\rho$ in the PR
representation \eqref{PRbiv} of the stochastic model \eqref{normal2},
which can never be identified from data. 

\section{Causal Inference from Experimental Studies}
\label{sec:expt}
Having set out a variety of formulations of my basic decision problem,
we now address the question: How might I gather and use data to help
identify the required ingredients?  From the DT perspective, I need to
assess my hypothetical distributions $P_0$ and $P_1$ for $Y$, under
either treatment choice.  These assessments should be informed by
(that is, conditioned on) whatever relevant information I may have:
for example, the responses observed for other similar headaches (my
own, or those of other people), that received one of the two
treatments.
  
We initially restrict attention to the simplest case.  Suppose I can
observe two groups of people.  Each group consists of individuals I
can regard as similar to (technically, exchangeable with) me in all
features relevant to their development of headaches and reaction to
treatment.  The first group are then all assigned active treatment
(aspirin): $X=1$; while the second group gets the control treatment
(no aspirin), $X=0$.\footnote{The usual operational method is, first
  to form a single group of individuals ``like me'', and, secondly, to
  randomise the assignment of treatment to its individuals; the
  resulting treatment/control groups will then each still be
  exchangeable with me on their pre-treatment characteristics.  The
  second stage by itself ensures {\em internal validity\/}: the
  treated and untreated groups should be comparable with each other;
  but without the first stage we will not have {\em external
    validity\/}, permitting generalisation beyond the data to external
  cases of interest --- in this case, myself.}  Finally I observe the
responses of all individuals: let $Y_{xi}$ denote the response of the
$i$th individual receiving treatment $x$.

\subsection{DT approach}
\label{sec:DTexp}
Under the above assumptions, I can model the responses of the treated
individuals as being randomly drawn from $P_1$, and likewise the
responses of the untreated individuals as drawn from $P_0$.  I can
then use completely standard statistical methods to estimate and
compare (in any way I choose) the two distributions $P_1$ and $P_0$.
In particular, I have access to all the ingredients required for my
stochastic decision problem.

For example, under model \eqref{normal2}, we model
\begin{equation}
  \label{eq:normdata}
  Y_{xi} \sim \norm(\mu_x, \sigma^2),
\end{equation}
all independently; and can then base inference about the difference
$\delta = \mu_1-\mu_0$ on Student's $t$-distribution.  All this is the
bread and butter of the most elementary courses in Statistics.  Here we
have however emphasised --- as may be more rarely done --- the
assumptions needed to justify the relevance of this inference to the
decision problem I face.

\subsection{Other approaches}
\label{sec:ESEexpt}
Suppose now we take the ESE model of \eqref{estruct} seriously (this
includes the simpler SSE model~\eqref{struct} as the special case $E_0
= E_1$, or, equivalently, $\rho=1$).  Because that model implies all
the distributional properties of \eqref{normal2}, I can still do all
that I did in \secref{DTexp} above, and so use the data to help solve
my decision problem.  But now I might want to do more.  For example, I
might want to say something about the distribution of my individual
causal effect, $\ice = \mu_1-\mu_0 + E_1 - E_0$.  Can I use the data
to help me in this?

Well, the mean of the $\ice$ is just $\delta = \mu_1-\mu_0$, which I
can estimate, as above.  But its variance is $2(1-\rho)\sigma^2$.
While I can estimate $\sigma^2$ from my data, the dependence on the
correlation $\rho$ is problematic: I could only estimate $\rho$ if I
had observations on bivariate pairs, $\bY = (Y_0,Y_1)$ ---
corresponding to observing both potential outcomes for the same
individual.  Since each individual receives just one of the two
treatments, the full pair $\bY$ is, logically, always unobservable,
and there is simply no way I can estimate $\rho$.  In particular, I
have no way of distinguishing observationally between the general ESE
model \eqref{estruct} and the SSE model~\eqref{struct}.  However,
these have different implications for $\var(\ice)$, which is
$2(1-\rho)\sigma^2$ for ESE, and $0$ for SSE.

In like fashion, were I to be interested in the mean of the ``ratio
\ice'', \viz\ $\E(Y_1/Y_0)$; or the estimate of my \ice after having
taken aspirin and observed response $Y=y$, \viz\ $\E(Y_1 - Y_0 \cd
X=1, Y = y) = (1-\rho)y + (\rho\mu_1-\mu_0)$; these seemingly
innocuous queries could not be addressed by any data I could ever
collect, since they all depend on the unknowable value of $\rho$.

As the ESE model is a special case of a PR model, or a general
functional model, all the above caveats apply to those models also.
We can thus divide putative causal inferences from such a model into
sheep --- those that depend only on the marginal distributions, $P_0$
and $P_1$, of the individual potential responses $Y_0$ and $Y_1$, and
are thus identifiable from data --- and goats --- those that do not,
and so are not identifiable.  For example, any putative causal
inference that makes essential use of $\var(\ice)$ is a goat.

It is all too easy to neglect this distinction, and attempt to make a
goat-like inference.  Given a fully specified ESE/PR/functional model
this will be mathematically possible, and it may not be noticed that
the answer would be different for a mathematically distinct model that
is observationally entirely equivalent to this one because it has the
same marginal distributions.  Such an apparent causal inference is, to
say the least, misleading.

\subsection{Treatment-unit additivity}
\label{sec:TUA}
The ESE model~\eqref{estruct} would represent the data as
\begin{equation}
  Y_{xi} = \mu_x + E_{xi}.
\end{equation}
For the special case of an SSE model~\eqref{struct}, this becomes
\begin{equation}
  \label{eq:estructdata}
  Y_{xi} = \mu_x + E_{i},
\end{equation}
a sum of one term, $\mu_x$, depending only on the treatment $x$
applied, and another, $E_i$, depending only on the unit $i$ to which
it is applied.  This property is termed {\em treatment-unit
  additivity\/} (TUA).  It is entirely equivalent to the property that
the $\ice$, $Y_1-Y_0$, has the identical value (\viz, $\mu_1-\mu_0)$
across all individuals.

One reason we might like the TUA assumption is as follows.  So far I
have had to assume that the individuals on whom I have data are ``like
me'' --- in particular, we all have the same distribution for our
error term $E$.  But this is often unrealistic, since \eg\ a clinical
trial will typically have stringent recruitment criteria, which I
would not satisfy.  However, I can relax model \eqref{estructdata}, so
as not to require that my own $E$ be drawn from the same distribution
as the $(E_i)$ in the data (for example, it could have a higher mean
or variance).  But since my own $\ice$ is still $\mu_1-\mu_0$, which
is still estimable from the data, the causal inference about my $\ice$
is unaffected by this relaxation of the SSE model --- though it does
still rely on TUA.

However an alternative decision-theoretic analysis (which thus does
not require TUA) is as follows.\footnote{See \textcite{apd:cinfer}
  \S~8.1 for a more detailed account.}  Because of different selection
criteria, my personal hypothetical distribution $P_x$ for my response
$Y$, if I take treatment $X=x$, is allowed be different from $P^*_x$,
the distribution of $Y$ for the individuals in my data who receive
treatment $X=x$.  But suppose (for example) I can model its mean
$\mu_x$ as related to the mean $\mu^*_x$ of $P_x^*$ by $\mu_x =
\mu^*_x +\gamma$ for some $\gamma$ that does not depend on $x$.  Then
my own \ace\ $\mu_1-\mu_0 = \mu^*_1-\mu^*_0$ is estimable from the
data.

\subsubsection{Neyman and Fisher}
\label{sec:neyfish}
Since there is no observational way of distinguishing between the SSE
model, for which TUA holds, and the more general ESE model, for which
it does not, the arguments in \secref{ESEexpt} would classify any
attempt at causal inference that is dependent on the assumption of TUA
as a goat.  In this light it is interesting to revisit the
ill-tempered debate between Neyman and Fisher in \textcite{jn:jrss}.

The essence of Neyman's model\footnote{We use different notation, and
  ignore certain elaborations irrelevant for current purposes.}
involves an experiment in which various treatments $t=1,\ldots,T$ can
be applied to various experimental units $u = 1,\ldots,U$, for example
plots in a field, which might themselves be described in more detail:
for example, in a ``randomised blocks'' layout, $u=(i,j)$
($i=1,\ldots.I; j=1,\ldots,J)$ for the $j$th plot in the $i$th row.
The treatments are applied to the units according to a randomisation
scheme taking account of the structure of the units.  In what is
considered to be the first use of a potential response formulation in
Statistics, Neyman introduces $y_{tu}$ to represent the response that
would be observed on unit $u$, if it were to receive treatment $t$ ---
with the values $(y_{tu})$ for fixed $u$, as $t$ varies, regarded as
having simultaneous existence, even though at most one can be
observed.  Neyman regards the collection of the $(y_{tu})$, for all
$t$ and $u$, as the ``unknown parameter''.  Statistical inference is
based on the distribution of the observed responses brought about by
the random assignment of treatments to units.

Neyman introduces the null hypothesis:
\begin{quote}
  $H^*_0:$ the value of $y_{t.}$ does not depend on $t$
\end{quote}
where $y_{t.}$ denotes the average of the $y_{tu}$ over the $U$ units
in the experiment.  That is, $y_{t.}$ is the average response that
would be obtained if treatment $t$ were applied to all the
experimental units, and his null hypothesis is that this would be the
same for every treatment --- allowing, however, that there could be
unit-level differences, that just happen to average out.\footnote{but
  need not do so if averaged over some other collection of units}
Neyman's analysis (corrected and extended by \textcite{w+k:jasa55})
shows that, for certain designs, such as the Latin square, the
standard $F$-test is a valid test of his $H^*_0$ only under the
assumption of treatment-unit additivity --- in which case, it may be
noted, $H^*_0$ becomes equivalent to
\begin{quote}
  $H^{**}_0:$ for each unit $u$, the value of $y_{tu}$ does not depend
  on $t$.
\end{quote}
Fisher criticised $H^*_0$, and the associated test, as based on an
entirely inappropriate formulation of the phrase ``no differences
between the treatments''.\footnote{Fisher's arguments are
  characteristically intuitive rather than formal --- though no less
  compelling for that --- and he is often taken as having favoured
  $H^{**}_0$, which is still phrased in terms of potential responses,
  as the appropriate null hypothesis.  However my own reading of his
  remarks does not find any clear commitment to a PR interpretation.}

From the point of view of the DT approach, it is troubling that,
according to Neyman, the standard $F$-test is or is not valid
according as whether or not we assume TUA --- a distinction without
any empirically observable consequences.  Neyman's analysis must
therefore be classified as a goat.

\section{Observational Studies and Confounding}
\label{sec:obsconf}
We have so far treated $X$ as a decision variable, under the control
of a human agent rather than Nature --- not merely for the actual
decision problem I myself face, but likewise for the experimental
individuals used to supply inputs for my problem.  But often genuine
experimentation is impossible, and we have to rely on data already
collected, in circumstances over which we had no control --- in
particular, no control as to who received which treatment.  Such {\em
  observational studies\/} raise serious problems of interpretation
and relevance, and great care is needed in drawing conclusions from
them \cite{rosenbaum:designbook,madigan:annrev}.

So suppose again I have data on a group of individuals whom I can
regard as exchangeable with me --- but now for whom treatments have
already been assigned, I know not how.  For each individual I have
information (say, for one headache episode) on the treatment applied
$X$, and the duration $Y$ (and, typically, also on some further
relevant variables).  Since I did not have the option to choose which
treatment to apply, $X$ is no longer a decision variable: it has
become a random variable.

A natural question is whether I can still use (an estimate of) the
distribution of $Y$, for those individuals who received treatment
$X=1$, as a proxy for my own hypothetical distribution $P_1$ (and
similarly for $X=0$).  Now in order for this to be possible, I must be
able to regard the {\em treated\/} patients as similar to
(exchangeable with) me, as regards relevant features existing prior to
treatment choice.  But even though this exchangeability may have been
assumed at the level of the whole group, it does not follow that it
will hold for the subgroup who got treated, since the treatment
decision may itself have been correlated with these features --- for
example, aspirin may have been given only for really bad headaches.
This is the problem of {\em confounding\/}, which obstructs
straightforward causal interpretation of observational data.  We shall
have {\em no confounding\/} only when I can, simultaneously, consider
myself as exchangeable, both with those patients who received aspirin,
and also with those who received none.  This then implies that those
two groups of patients must be exchangeable with each
other,\footnote{That is, external validity for each group separately
  implies internal validity.} which in turn requires that the way in
which treatments were applied was oblivious to (independent of) any
features of the individuals that could be relevant to their reactions
to treatment.  The easiest way to ensure this is by randomisation.
Although we are here assuming this was not possible, we can sometimes
(albeit rarely) attempt to argue that the data can nevertheless be
treated as if they had been randomised.

For a functional model $Y=f(X,U)$, the defining property of ``no
confounding'' is typically taken as requiring independence (in the
observational data) between $X$ and $U$: in the notation of
\textcite{apd:CIST} (see \secref{ci} below),
\begin{equation}
  \label{eq:xindu}
  \indo X U.  
\end{equation}
This is trivially equivalent to $\indo U X$, \ie\ the observational
distribution of $U$ given $X=x$ does not depend on $x$, thus mimicking
a property already assumed for the case that $X$ is my decision
variable.  For an ESE model the above requirement translates as $\indo
X \bE$, and as $\indo X \bY$ for a PR model.  (However, since $U$,
$\bE$ and $\bY$ typically do not correspond to any empirically
observable variables, the mental exercise required to assess whether
the above independence properties hold can be perplexing.)

In the next section we consider just why \eqref{xindu} might be
considered as expressing ``no confounding'', and extend the analysis
to the DT interpretation of this concept.

\subsection{Regimes}
\label{sec:regimes}

A helpful way to think about (absence of) confounding is in terms of
different data-generating {\em regimes\/} and their relationships.

In the above example we can distinguish three such regimes.  One of
these is the {\em observational\/} regime, under which the available
data have been observed.  Then there are two {\em interventional\/ }
regimes, corresponding to the circumstance in which an external
intervention is made to impose one of the two treatments.  It is
helpful to introduce a non-stochastic variable $F_X$, with values $0$,
$1$ and \idle\ (read as ``idle''): $F_X=0$ labels the interventional
regime with $X=0$; $F_X=1$ labels the interventional regime with
$X=1$; while $F_X = \idle$ labels the observational regime.  There
will be a joint distribution of all relevant variables for each of
these regimes.  Thus $F_X$ has the status of a statistical parameter,
indexing which distribution we are referring to.  Note that (assuming
intervention is perfectly successful), $X=x$ with probability 1 in
regime $F_X = x$ ($x = 0, 1)$, whereas $X$ will be a genuinely
stochastic variable in regime $F_X = \idle$.

We have previously interpreted a functional model, and its
specialisations such as an ESE or PR model, as incorporating an
implicit assumption of ``stability'': that the relevant variable $U$
should have the same value, and hence the same distribution, no matter
which treatment is applied.  In fact this assumption is not quite
enough: we further need to assume that $U$ has the same distribution,
irrespective of the regime that is operating.  This is necessary if we
are to justify transfer (under suitable conditions) of information
from the observational regime to the interventional regimes.  Thus
suppose \eqref{xindu} holds.  We desire to compute, and contrast, the
distributions of the response $Y$ under the two interventional
regimes.  Under intervention with active treatment, $Y = f(1,U)$ with
$U$ having its distribution under $F_X = 1$.  In the observational
regime, we can estimate the conditional distribution of $Y$ given
$X=1$, which is that of $f(1,U)$ given $X=1$ under $F_X = \idle$.  On
account of \eqref{xindu} (supposed to apply in the observational
regime), this is the same as the marginal distribution of $f(1,U)$
under $F_X = \idle$.  But under our extended stability assumptions,
$U$ has the same distribution in all regimes, so this is indeed the
same as the desired distribution of $Y = f(1,U)$ under $F_X = 1$.  We
have thus shown that, taking together \eqref{xindu} and the
assumptions of stability, we can deduce ``no confounding'', expressed
as
\begin{equation}
  \label{eq:sim}
  (Y \cd F_X = x) \approx (Y \cd X=x; F_X = \idle) \quad(x=0,1),
\end{equation}
where the symbol $\approx$ denotes ``has the same distribution as''.

At this point we note that the left-hand side of \eqref{sim} refers to
what we have termed the ``hypothetical distribution'', $P_x$, of $Y$,
under intervention with $X=x$; while the right-hand side refers to a
conditional distribution that is, in principle, estimable from
empirical data.  All the special ingredients of the functional model
have evaporated, and we are left with an expression that is fully
meaningful within the DT framework.  And within that framework we can
simply and directly take property \eqref{sim} (however it may be
justified) as the appropriate expression of ``no confounding''.

\subsection{Conditional independence}
\label{sec:ci}
An alternative way of expressing \eqref{sim} is as follows.  First
note that, since $X=x$ with probability 1 under regime $F_X = x$,
\eqref{sim} is is equivalent to
\begin{equation}
  \label{eq:sim1}
  (Y \cd X=x; F_X = x) \approx (Y \cd X=x; F_X = \idle) \quad(x=0,1).
\end{equation}
This expresses the conditional distribution of $Y$, given $X=x$, as a
``modular component'', that can be transferred without change betwen
observational and interventional settings.  This modular
interpretation of ``causality'' offers a useful pragmatic take on a
slippery philosophical concept.

The distributional identity \eqref{sim1} can also be considered as an
expression of the {\em conditional independence\/} property
\cite{apd:CIST,apd:ciso}:
\begin{equation}
  \label{eq:ciprop}
  \ind Y {F_X} X,
\end{equation}
which says that the distribution of $Y$, given information both on the
value of $X$ and on the regime $F_X$ under which that value arose, is
in fact the same for all the regimes.  In this way we have converted a
causal property into a probabilistic one (albeit involving the
non-random regime variable $F_X$).  Since there is a well-established
theory of conditional independence (see \secref{ci} below), this is a
fruitful reinterpretation that will be particularly helpful for both
describing and manipulating causal properties.  So henceforth we will
work with \eqref{ciprop}, and its DT interpretation, as our formal
expression of ``no confounding''.

\section{Conditional Independence and Graphs}
\label{sec:ci}

In this section we recapitulate various aspects of the mathematical
theory of conditional independence that will be useful for
manipulating causal concepts.  For further detail see
\textcite{apd:CIST,apd:ciso,apd:infdiags,nayia:thesis}.

For random variables $X$, $Y$, \ldots, having joint distribution $P$,
we say {\em $X$ is independent of $Y$ given $Z$\/}, and write $\ind X
Y Z$, to mean that the distribution of $X$, given $(Y,Z)=(y,z)$
depends only on the value $z$ of $Z$.  More formally:
\begin{defn}[{\sc Conditional independence}]
  \label{def:CI-ci}
  We say {\em $X$ is conditionally independent of $Y$ given $Z$\/},
  and write $\ind X Y Z$, if, for any measurable set $A$ in the range
  of $X$, there exists a function $w(Z)$ of $Z$ alone such that $P(X
  \in A \cd Y, Z) = w(Z)$\ [$P$-almost surely].
\end{defn}
When we need to specify explicitly the underlying joint distribution
$P$ we write \eg\ $\indp X Y Z P$.  Independence, $\indo X Y$, is the
special case of conditional independence when the conditioning
variable $Z$ is trivial.

\subsection{Axioms of conditional independence}
\label{sec:CI-prax}

Among the general properties of probabilistic conditional independence
are the following \cite{apd:CIST}.  Here we write $W\preceq Y$ to mean
that $W$ is a function of $Y$.

\begin{displaymath}
\begin{array}{rlclcl}
\mbox{P1} & \mbox{\em ``Symmetry''} & : & \ind X Y Z &  \Rightarrow &
\ind Y X Z \\
\mbox{P2} && : & \ind X Y X  \\
\mbox{P3} & \mbox{\em ``Decomposition''}& : & \ind X Y Z,\quad  W\preceq
Y & \Rightarrow & \ind X W Z \\
\mbox{P4} & \mbox{\em ``Weak union''} & : & \ind X Y Z,\quad  W\preceq Y
& \Rightarrow & \ind X  Y {(W,Z)} \\
\mbox{P5} & \mbox{\em ``Contraction''}& : & \!\!\!\left. \begin{array}{l}
   \ind X Y Z \\
   \mbox{\hspace{1em}and} \\
   \ind X W {(Y,Z)}
   \end{array}
\right\} &  \Rightarrow & \ind X {(Y,W)} Z.
\end{array}
\end{displaymath}
(The descriptive terms are those given by
\textcite{pearl:pesbook},~Chapter 3).

It is possible to derive many further properties of CI by regarding P1
to P5 as axioms for a logical system, rather than calling on more
specific properties of probability distributions.

\subsection{Extension to non-stochastic variables}
\label{sec:CI-nonstoch}
In order for \defref{CI-ci} to make sense we must be able to talk
about distributions for $X$, which thus has to be a random variable;
but (subject to appropriate interpretation of the ``almost sure''
qualification) $Y$ and $Z$ need not be.  In particular, this is the
case for expression \eqref{ciprop}, our interpretation of ``no
confounding'', which involves the non-stochastic regime indicator
variable $F_X$.

We must exercise a little care when applying the notation and theory
of \secref{CI-prax} to non-stochastic variables, to ensure that these
always appear, explicitly or implicitly, as conditioning variables.
Nevertheless, suitably interpreted, properties P1--P5 do still hold
\cite{apd:ciso,nayia:thesis}.  In fact any deduction made using them
will be valid, so long as, in both premisses and conclusions, no
non-stochastic variables appear in the left-most term in a conditional
independence statement (we {\em are\/} allowed to violate this
condition in intermediate steps of an argument).  So we can apply
P1--P5 freely, even in the presence of non-stochastic variables, so
long only as we do not attempt to derive any obviously meaningless
assertion.

\subsection{Graphical representation}
\label{sec:graphs}

There is a remarkable and technically valuable analogy between
conditional independence properties holding between random variables,
and separation properties of a directed acyclic graph (DAG)
\cite{sll/apd/bnl/hgl:directed}.  This enables us to use graphical
methods to streamline probabilistic manipulations.

The graphical analogue of probabilistic conditional independence is
the following somewhat complex separation property.  Let $A$, $B$, $C$
be sets of nodes of the DAG ${\cal D}$.  We first form the subgraph
$\dag'$ of $\dag$ that contains only the nodes in $A$, $B$ and $C$,
together with all their ancestors in $\dag$, and all their connecting
arrows: this is the relevant {\em ancestral DAG\/}.  Next, whenever
two nodes in $\dag'$ have a common child but are not already joined by
an arrow (are ``unmarried''), we insert an undirected edge between
them, and then convert all remaining edges to be undirected by
dropping the arrowheads: this produces the {\em moralised\/} ancestral
graph, ${\cal G}'$.  Finally, in the undirected graph ${\cal G}'$, we
check whether every connected path from a node in $A$ to one in $B$
intersects $C$.  If so, we say {\em $C$ $d$-separates $A$ from
  $B$\/},\footnote{The name refers to an alternative, but equivalent,
  way of expressing this separation property, as described by
  \textcite{pearl:blocking,verma}.} and write $\inda A B C \dag$.  It
can be shown that, at a purely formal level, and with $\preceq$ now
interpreted as ``is a subset of'', $\subseteq$, this separation
property satisfies Axioms P1--P5 of \secref{CI-prax}

Now suppose that each node $v$ of $\dag$ has an associated random
variable $X_v$.  Denote $(X_v: v\in A)$ by $X_A$.  We say a joint
distribution $P$ for all these variables satisfies the {\em local
  directed Markov property\/} with respect to $\dag$ if, for every
node $v$,
\begin{equation}
  \label{eq:DAG-parci}
  \ind {X_v}  {X_{\nd{v}}} {X_{\parents{v}}},
\end{equation}
where (using a self-explanatory analogy with a genetic pedigree)
$\parents{v}$ denotes the set of ``parents'' of node $v$, and $\nd{v}$
its ``non-descendents'', in $\dag$.  In this case it can be shown
that, whenever we find $\inda A B C \dag$ (by inspection of the DAG),
we can deduce the probabilistic conditional independence property
$\indp {X_A} {X_B} {X_C} P$.  We term this the {\em moralisation
  criterion\/}.

\begin{figure}
\begin{center}
  \resizebox{3in}{!}{\includegraphics{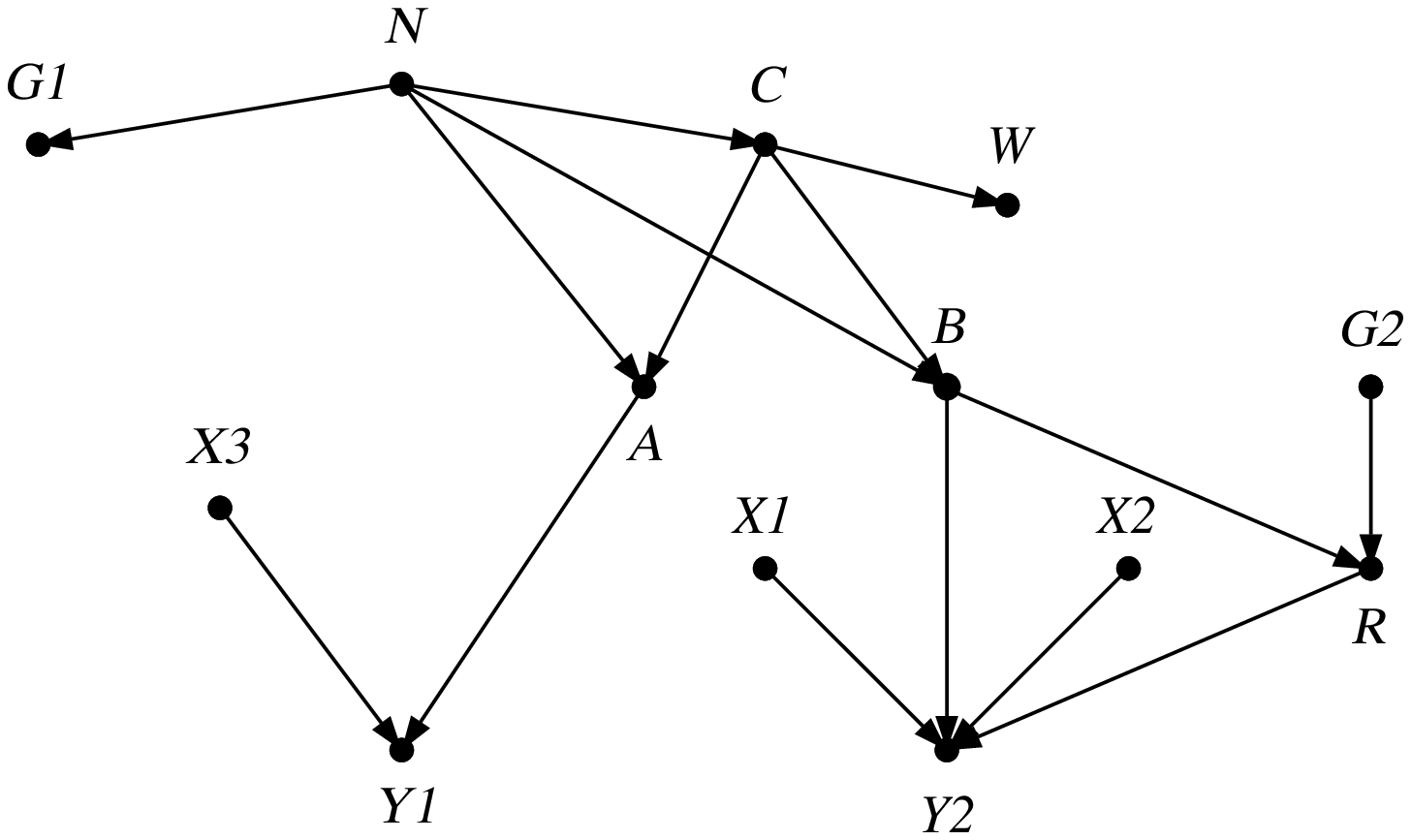}}
\caption{Directed graph $\cal D$ for criminal evidence}
\label{fig:dag}
\end{center}
\end{figure}

As an example, the directed acyclic graph (DAG) $\cal D$ of
Figure~\ref{fig:dag} describes the relationships between the evidence
and other variables figuring in a criminal trial \cite{apd/iwe}.
The graph is constructed so that the each node corresponds to a
variable in the problem, and the assumed dependence structure of the
variables satisfies the local directed Markov property: each variable
is supposed probabilistically conditionally independent of its
non-descendents in the graph, conditional on its graph parents.  For
example, the distribution of $Y1$ (measured properties of a tuft of
fibres found at the scene), given all other variables, is supposed
fully determined by the values of $X3$ (properties of the suspect's
jumper) and of $A$ (an indicator of whether or not the fibres came
from the suspect's jumper).  Likewise, the distribution of $B$ (who
left blood on the jumper?), given all variables other than $Y2$ (the
type of that blood) and $R$ (whether or not the blood pattern was a
spray), in fact only depends on the values of $N$ (the number of
offenders) and $C$ (whether or not the suspect was an offender).  Such
assessments can often be made at a qualitative level, before
attempting numerical specification of probabilities.  In turn, that
specification is simplified because we only need to describe the
conditional distribution for each variable given its graph parents.
 
\begin{figure}
  \begin{center}
    \resizebox{3in}{!}{\includegraphics{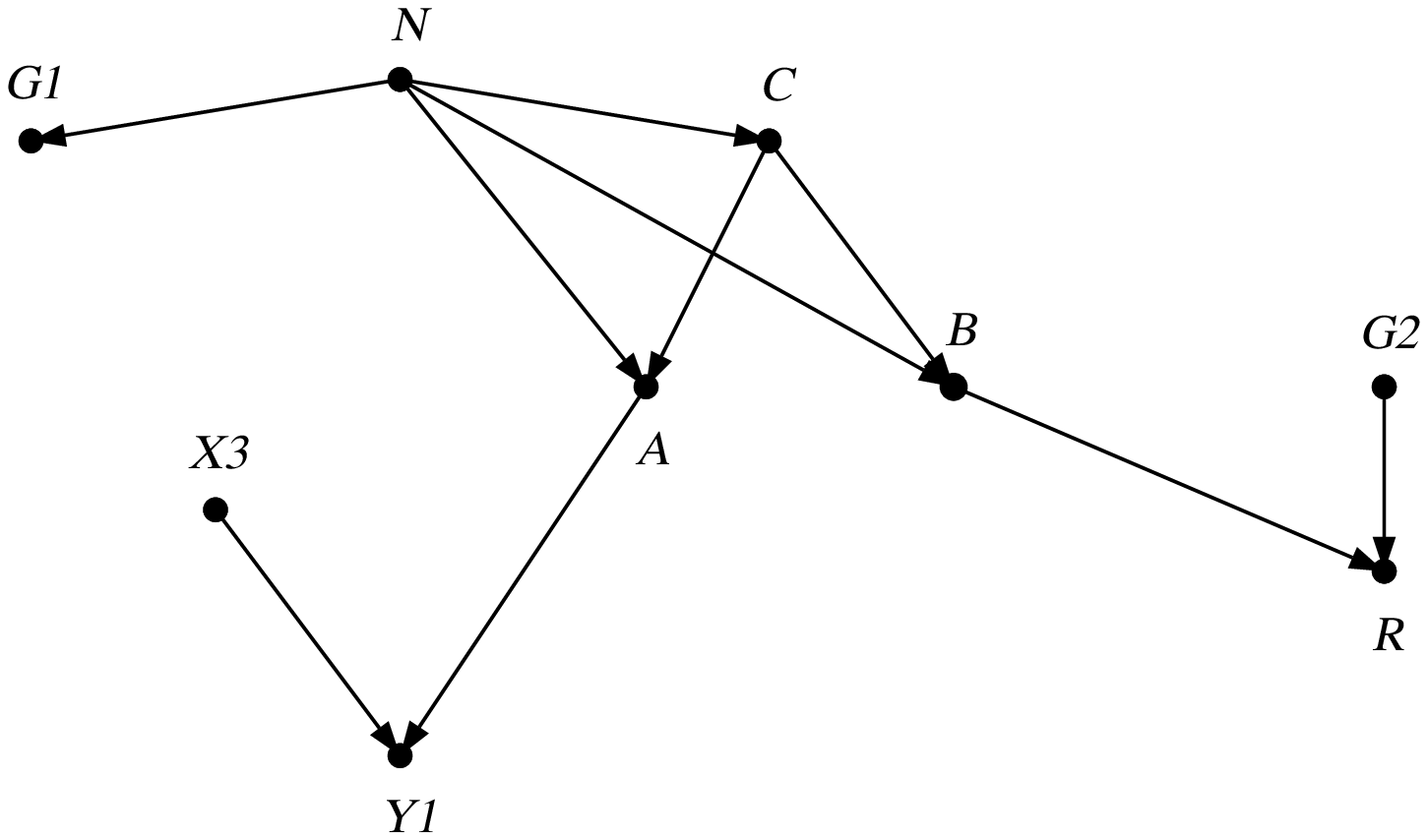}}\end{center}
\caption{Ancestral subgraph ${\cal D}'$}
\label{fig:anc}
\end{figure}

\begin{figure}
\begin{center} \resizebox{3in}{!}{\includegraphics{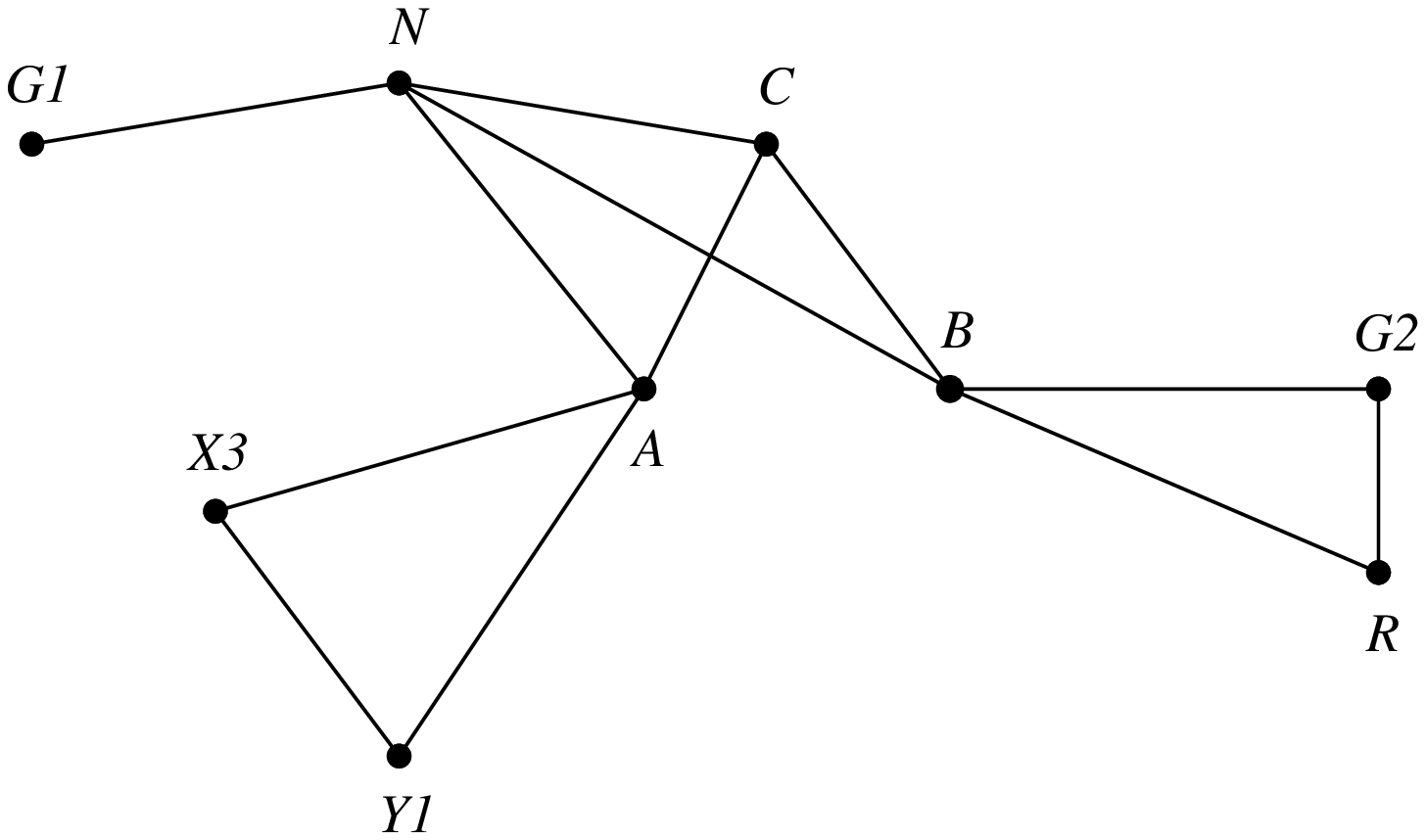}}\end{center}
\caption{Moralised ancestral subgraph ${\cal G}'$}
\label{fig:mor}
\end{figure}

Suppose now we wish to query whether $\ind {(B,R)} {(G1,Y1)} {(A,N)}$.
The relevant ancestral graph $\dag'$ is shown in \figref{anc},
and its moralised version ${\cal G}'$ in \figref{mor}.
We note that in ${\cal G}'$ it is impossible to trace a path from
either of $B$ or $R$ to either $G1$ or $Y1$ without passing through
either $A$ or $N$.  Thus $\inda {(B,R)}{(G1,Y1)}{(A,N)} \dag$.  From
this we deduce the probabilistic conditional independence property
$\ind {(B,R)} {(G1,Y1)} {(A,N)}$.

\paragraph{Caution:} Although every DAG thus describes some collection
of conditional independence properties, and can be used to manipulate
these, by no means every such collection can be represented by a DAG.
In full generality, we may need to use algebraic manipulations,
successively applying the CI axioms P1--P5 to derive the implicit
consequences of an assumed collection of conditional independencies.

\subsubsection{Markov equivalence}
\label{sec:DAG-markeq}

Distinct DAGs can have identical separation properties, and so
represent identical collections of conditional independencies.  They
are then termed {\em Markov equivalent\/}.

The {\em skeleton\/} of a DAG $\dag$ is the undirected graph obtained
by ignoring the directions of the arrows on the edges of $\dag$.  An
{\em immorality\/} in $\dag$ is a configuration of the form $a \to c
\leftarrow b$, where $a$ and $b$ are parents of a common child $c$ but
neither $a \to b$ nor $b \to a$.

\begin{theorem}[\textcite{morten:markov,verma/pearl:90}]
  \label{thm:DAG-markoveq}
  Two DAGs $\dag_0$ and $\dag_1$ on the same vertex set $V$ are Markov
  equivalent if and only if they have the same skeleton and the same
  immoralities.
\end{theorem}

\begin{ex}
  \label{ex:DAG-markeq0}
  There are just three possible DAGs on two nodes:
  \begin{enumerate}
  \item \label{it:DAG-markeq02} $A\to B$
  \item \label{it:DAG-markeq03} $A\leftarrow B$
  \item \label{it:DAG-markeq01} $A\mbox{\hspace{3.5ex}}B$.
  \end{enumerate}
  Since DAGs \itref{DAG-markeq02} and \itref{DAG-markeq03} have the
  same skeleton, and neither has any immoralities, they are Markov
  equivalent: indeed, they embody no conditional independence
  properties whatsoever.  However, DAG \itref{DAG-markeq01}, which has
  a different skeleton, embodies the non-trivial conditional
  independence restriction $\indo A B$.
  \end{ex}
  
\begin{ex}
  \label{ex:DAG-markeqsimp}
  Consider the following DAGs on three nodes:
  \begin{enumerate}
\item 
      \label{it:DAG-markeqsimp1}
      $A \to B \to C$
\item 
      \label{it:DAG-markeqsimp2}
      $A \leftarrow B \leftarrow C$
\item 
      \label{it:DAG-markeqsimp3}
      $A \leftarrow B \to C$
\item 
      \label{it:DAG-markeqsimp4}
      $A \to B \leftarrow C$.
  \end{enumerate}
  These all have the same skeleton.  However, whereas DAGs
  \itref{DAG-markeqsimp1}, \itref{DAG-markeqsimp2} and
  \itref{DAG-markeqsimp3} have no immoralities,
  \itref{DAG-markeqsimp4} has one immorality.  Consequently,
  \itref{DAG-markeqsimp1}, \itref{DAG-markeqsimp2} and
  \itref{DAG-markeqsimp3} are Markov equivalent to each other, but
  \itref{DAG-markeqsimp4} is not Markov equivalent to these.  Indeed,
  \itref{DAG-markeqsimp1}, \itref{DAG-markeqsimp2} and
  \itref{DAG-markeqsimp3} all express the conditional independence
  property $\ind A C B$, whereas \itref{DAG-markeqsimp4} expresses the
  marginal independence property $\indo A C$.
\end{ex}

\section{Causal Interpretations of DAGs}
\label{sec:CAUSDAG}

\begin{figure}[htbp]
  \begin{center}
    \resizebox{1.5in}{!}{\includegraphics{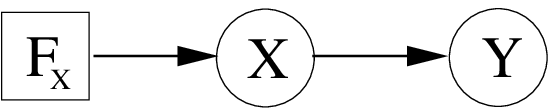}}
   \caption{No confounding: $X$ causes $Y$}
   \label{fig:x2y}
  \end{center}
\end{figure}

It is common, and appears very natural, to want to interpret an arrow
$a \to b$ in a DAG as representing some kind of ``direct causal
dependence'' of $b$ on $a$.  But this is a potentially dangerous move,
since there is nothing in the DAG semantics, as presented above, to
justify it.  We prefer a different way of introducing causality into a
DAG: by explicitly representing regime indicators,\footnote{This
  develops on an idea introduced by \textcite{sgs:book}: see also
  \textcite{pearl:book,sll:causal}.} and applying the moralisation
criterion to the resulting {\em influence diagram\/} (ID), a DAG
containing both stochastic and non-stochastic variables.  As a simple
example, the ``no confounding'' property \eqref{ciprop} is represented
by the ID of \figref{x2y}.

\begin{figure}[htbp]
  \begin{center}
    \resizebox{1.5in}{!}{\includegraphics{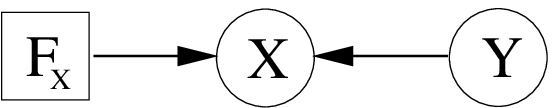}}
    \caption{$X$ does not cause $Y$}
    \label{fig:y2x}
  \end{center}
\end{figure}

Consider now the effect of reversing the arrow from $X$ to $Y$, as
shown in \figref{y2x}.
Without the intervention node $F_X$, the two graphs would have been
Markov equivalent (as was the case for \exref{DAG-markeq0},
\itref{DAG-markeq02} and \itref{DAG-markeq03}).  Now however we can
easily see that they no longer represent equivalent assumptions since,
although they have the same skeleton, they have different
immoralities.  \figref{y2x} expresses the marginal independence
property $\indo Y {F_X}$, and thus makes it explicit that the {\em
  marginal\/} distribution of $Y$ is the same, no matter whether, or
how, $X$ is subjected to intervention.  That is, $X$ has no effect on
$Y$ (in any regime).

\section{Pearlian DAGs}
\label{sec:pearlian}

\begin{figure}[htbp]
  \begin{center}
    \resizebox{1.8in}{!}{\includegraphics{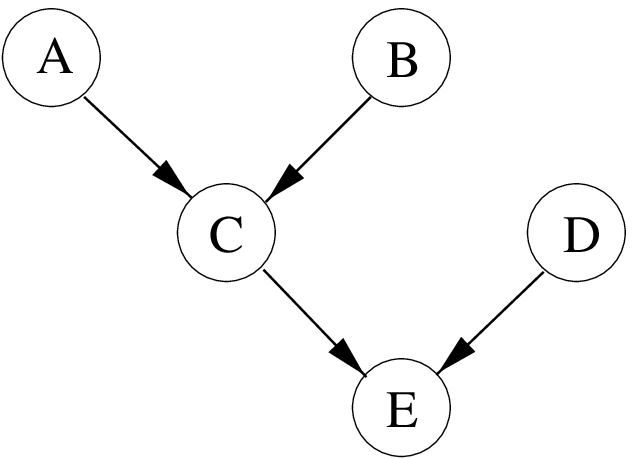}}
    \caption{A  probabilistic DAG}
    \label{fig:probdag}
  \end{center}
\end{figure}

Consider the DAG of \figref{probdag}.
Interpreted purely stochastically, it is nothing but a representation
of the following conditional independence properties: $\indo A B$;
$\indo D {(A,B,C)}$; and $\ind E {(A,B)} {(C, D)}$; together with all
other properties, such as $\ind E B {(A,C)}$, deducible from these
using P1--P5 (or, equivalently, readable off the DAG using the
moralisation criterion).

\begin{figure}[htbp]
  \begin{center}
    \resizebox{3.2in}{!}{\includegraphics{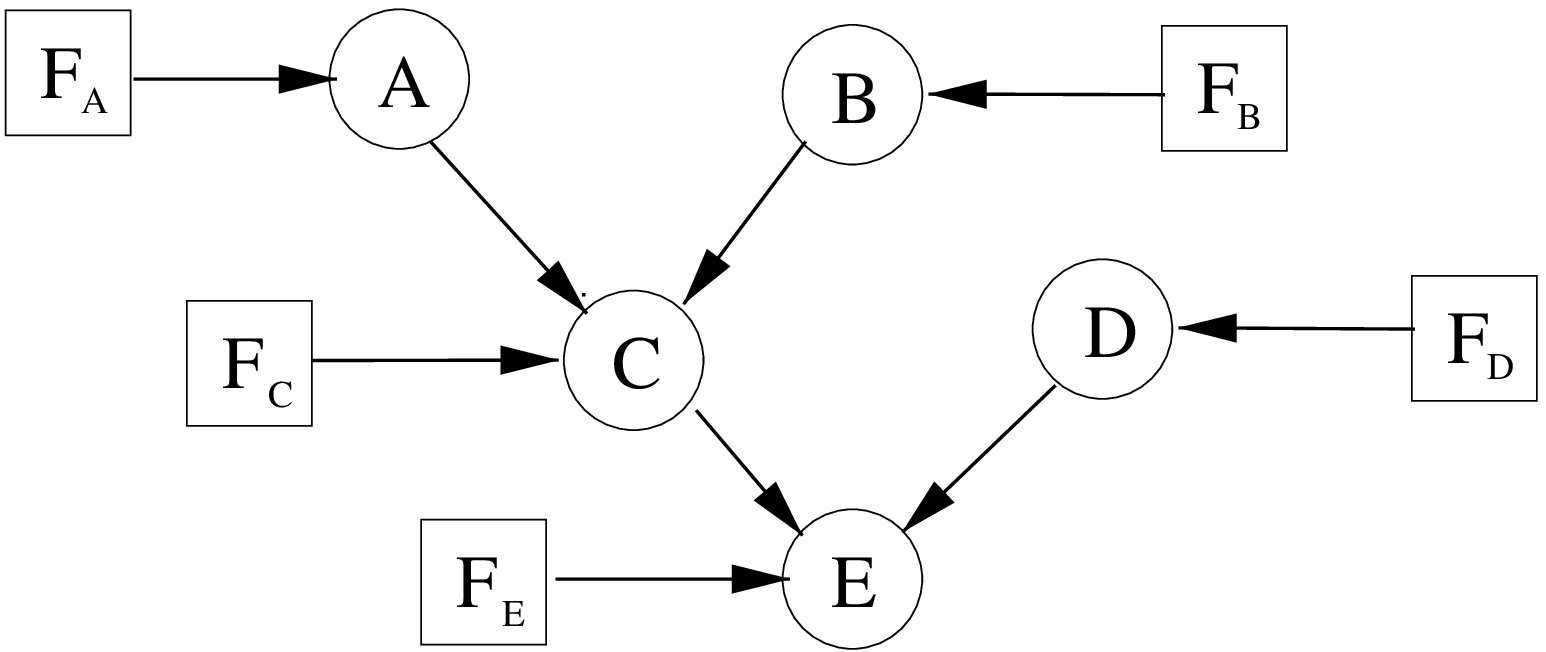}}
    \caption{Augmented DAG}
    \label{fig:augdag}
  \end{center}
\end{figure}

In the approach of \textcite{pearl:book}, a DAG such as
\figref{probdag} is taken to represent causal properties.  A helpful
way of understanding Pearl's interpretation is to consider the DAG as
a shorthand for the influence diagram of \figref{augdag}, in which a
non-stochastic intervention node has been associated with every
stochastic node.
Using the moralisation criterion, we can read off from this {\em
  augmented DAG\/} that, for example, $\ind C {(D, F_A, F_B, F_D,
  F_E)} {(A, B, F_C)}$.  For $F_C = \idle$ (the only non-trivial
case), this says that the `natural' conditional distribution of $C$,
given $A$ and $B$, is not further affected by additional conditioning
on the value of $D$, {\em nor by whether or not any or all of $A$,
  $B$, $D$ or $E$ arose naturally or by intervention\/}.  Similar
properties hold for the other domain variables.  In particular, we can
see that the conditional distribution for a node, given its domain
parents, when it is allowed to arise naturally, remains unchanged when
its parents are set by intervention (and is thus a modular component,
invariant across different regimes).  The augmented DAG thus
automatically encodes ({\em via\/} moralisation semantics) the
assumptions made externally by Pearl, without requiring any new
ingredients or concepts; and further makes it easy to read off their
implications directly.  It also makes it clear that, when endowed with
Pearl's causal interpretation, DAGs that are {\em prima facie\/}
Markov equivalent (such as $X\rightarrow Y$ and $X\leftarrow Y$) are
not causally equivalent, since their augmented forms will not be
Markov equivalent.  For all these reasons, explicit use of augmented
DAGs is to be preferred over Pearl's shorthand form, which in any case
courts confusion with the purely stochastic interpretation of a DAG.
 
\paragraph{Caution:}
A Pearlian DAG model, or its augmented DAG equivalent, is justified
only to the extent that it models the actual the behaviour of the
world in the setting to which it is intended to apply.  In particular,
we must ask whether or not the various interventional situations are
indeed related to the non-interventional one in the specific way
represented by the DAG.  Since such considerations necessarily involve
cross-regime comparisons, no assessment of their appropriateness can
be made on the basis of purely observational data.  

\section{Identifying Causal Effects}
\label{sec:ident}

Suppose we are interested in the ``causal effect'' of a treatment
variable $T$ on a response variable $Y$.  In the DT framework, this
requires us to identify, and contrast, the two interventional
distributions: $P_1$, for $Y$ in regime $F_T=1$, and $P_0$, for $Y$ in
regime $F_T=0$.  For simplicity we again confine attention to the
average causal effect
\begin{equation}
  \label{eq:ace}
  \ace: = \E(Y \cd F_T=1) - \E(Y \cd F_T=0).
\end{equation}
With only observational data, gathered in regime $F_T = \idle$, we
will not be in a position directly to assess these interventional
distributions of $Y$.  We will thus need to make assumptions to
justify and guide computation of the \ace from such data.  Since any
such assumptions will have to relate distributions across distinct
regimes, they will not be empirically testable if we only have
observational data.  It will however be important to present some sort
of convincing argument for the suitability of any assumptions imposed.

At the simplest level, we might assume ``no confounding'': $\ind Y
{F_T} T$.  In this case we could simply estimate the observational
conditional distribution of $Y \cd T=t, F_T=\idle$, and take that as
the desired interventional distribution of $Y \cd F_T=t$ ($t=0,1)$.
Thus under this assumption we will have
\begin{equation}
  \label{eq:acenoconf}
  \ace = \E(Y \cd T=1, F_T=\idle) - \E(Y \cd T=1, F_T=\idle),
\end{equation}
which is straighforwardly estimable from observational data.

However, in many realistic contexts the ``no confounding'' property
will be simply unbelievable: that is to say, we will have {\em
  confounding\/}: $\nind Y {F_T} T$.  Then \eqref{acenoconf} might
fail.  Note that this definition of confounding does not require the
existence of what are often called {\em confounding variables\/}, or
{\em confounders\/}.  But to make progress in identifying \ace\ we
will typically have to introduce further variables, with appropriate
properties.

\subsection{Sufficient covariates}
\label{sec:suffcov}
A variable $U$ is a {\em pre-treatment variable\/} if it exists and is
(in principle) observable prior to the point at which the treatment
decision is taken.  In this case its value, and so its distribution,
must be the same under both interventional regimes, $F_T = 0$ and $F_T
= 1$.  It will frequently (though not invariably) also be the case
that its distribution can be considered the same in the relevant
observational regime $F_T = \idle$.  Then we will have
\begin{equation}
  \label{eq:prop1}
  \indo U { F_T}.
\end{equation}
Such a variable is a {\em covariate\/}.

\subsection{Unconfounder}
\label{sec:unconf}
When we can not assume ``no confounding'', we might be able to tell an
alternative, more convincing, story, in terms of a (typically
multivariate) covariate $U$: claiming that we will have {\em no
  residual confounding\/}, after conditioning on $U$.  Formally,
\begin{equation}
  \label{eq:prop2}
  \ind Y {F_T} {(U,T)}.
\end{equation}
For example, if our data arise from an observational study on patients
treated by a certain doctor, who might be allocating treatment
according to his own observations $U$ of the general health of the
patient, it could be reasonable to suppose that, conditionally on $U$,
we would have no residual confounding.  If we can observe $U$, we can
then use the observational distribution of $Y$ given $(U, T=t)$ as the
distribution of $Y$ given $U$ in the interventional regime $F_T = t$.

A variable satisfying both \eqref{prop1} and \eqref{prop2} is often
called a {\em confounder\/}, though a more appropriate term might be
{\em unconfounder\/}.  We shall call it a {\em sufficient
  covariate\/}.

\begin{figure}[htbp]
  \begin{center}
    \resizebox{1.6in}{!}{\includegraphics{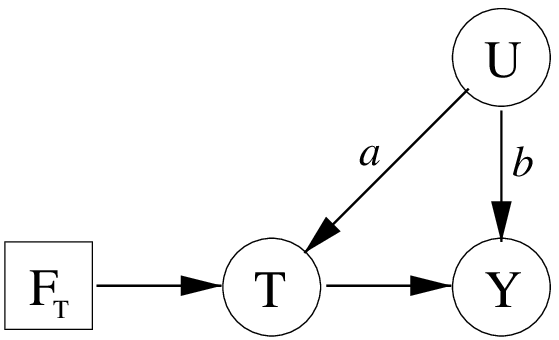}}
    \caption{Sufficient covariate}
    \label{fig:confsimpab}
  \end{center}
\end{figure}

The properties \eqref{prop1} and \eqref{prop2} are represented by the
ID of \figref{confsimpab}.

\subsubsection{Functional model}
\label{sec:fmodcon}
Suppose our starting point was a functional model $Y=f(T,U)$ (which
includes (E)SE and PR models).  Since the same function is supposed to
apply irrespective of the regime operating, \eqref{prop2} holds
trivially.  We have further assumed that $U$ has the same value, and
hence the same distribution, in all regimes, so \eqref{prop1} holds.
That is, formally at least, $U$ is an unconfounder.  However, in such
a formulation the variable $U$ is typically unobservable (this being a
logical necessity in the PR approach, where $U$ is the pair $\bY$ of
potential responses), which limits the operational usefulness of this
observation.

\subsection{Non-confounding}
\label{sec:nonconf}

Specialisations of the above structure are obtained when we can assume
that either of the arrows marked $a$ and $b$ in \figref{confsimpab} is
absent.  It can readily be checked that in either case we will have
$\ind Y {F_T} T$: no confounding.  We might call such a sufficient
covariate $U$ a {\em non-confounder\/}, and can safely forget that it
ever existed: we can simply apply \eqref{acenoconf}.

The ID with arrow $a$ absent represents the additional\footnote{\ie,
  over and above properties \eqref{prop2} and \eqref{prop1}.}
property $\ind T U {F_T}$~: that in every regime $T$ is independent of
$U$.  Since this condition holds trivially for the interventional
regimes, where $T$ is constant, it merely requires that $T$ be
independent of $U$ in the observational regime --- that is, that the
variables $U$ that putatively might have affected the doctor's
decision did not in fact do so.  This property would be perfectly
believable if the doctor had tossed a coin to determine his decision,
which is why randomised studies can directly address causal queries.
But in the case of an observational study, we would need to make some
alternative convincing case for this property: then (and only then) we
can treat the study as if it had been randomised.  This argument is
similar to that of \secref{regimes} for functional models, but ---
since it involves a real rather than a fictitious variable $U$, and
stochastic rather than deterministic relationships --- supplies a more
operational justification for assuming ``no confounding''.

The ID with $b$ absent represents the additional property $\ind Y U
T$, which says that the conditional distribution of $Y$ given $(T, U)$
(which, by \eqref{prop2}, has already been supposed the same in all
regimes) does not in fact depend on $U$: that is, $U$ is not
predictive of outcome.  In that case, even if $U$ is associated with
treatment assignment, this will not generate confounding.

\subsection{Deconfounding}
\label{sec:deconf}

More generally, suppose $U$ is a sufficient covariate that is observed
in the observational regime.  Define
\begin{eqnarray}
  \label{eq:sce}
  \sce_U &:=& \E( Y \cd U, F_T=1) - \E( Y \cd U, F_T=0),
\end{eqnarray}
the {\em specific causal effect\/} of treatment, given $U$.  This is a
random variable, a function of $U$, whose value $\sce_U(u)$ when $U=u$
is the average treatment effect in the subgroup of individuals having
$U=u$.

Now $T=t$ with probability 1 under $F_t$.  Then using \eqref{prop2} we
find $\E(Y \cd U, F_T=t) = \E(Y \cd U, T=t, F_T=t) = \E(Y \cd U, T=t,
F_T=\idle)$.\footnote{More accurately, these identifications require
  an additional {\em positivity\/} condition
  \cite{hg/apd:aistats2010}, which will typically be satisfied.}  We
deduce
\begin{equation}
  \label{eq:sceobs}
  \sce_U = \E(Y \cd  U, T=1, F_T=\idle)-\E(Y \cd  U, T=0, F_T=\idle),
\end{equation}
so that $\sce_U$ is estimable from observational data.  This is a
reflection of the fact that we have no confounding conditional on $U$.

Also, by the ``extension of the conversation'' rule of probability, we
have
\begin{eqnarray*}
  \E(Y \cd F_T=t) &=& \E\{\E(Y \cd  U, F_T=t) \cd  F_T=t\}\\
  &=& \E\{\E(Y \cd  U, F_T=t) \cd  F_T=\idle\}
\end{eqnarray*}
by \eqref{prop1}.  It follows that
\begin{equation}
  \label{eq:aceice}
  \ace = \E(\sce_U \cd F_T=\idle).
\end{equation}
That is, for any sufficient covariate $U$, the overall average causal
effect is the observational expectation of the associated specific
causal effect.  Since, by \eqref{sceobs}, $\sce_U$ is itself an
observationally estimable quantity, formula~\eqref{aceice} allows us
to estimate $\ace$ whenever we can observe a sufficient covariate.

Note that in the PR framework, where we take $U = \bY$, \sce\ becomes
$Y_1 - Y_0$, the ``individual causal effect'', \ice.  Then
\eqref{aceice} shows that $\ace = \E(\ice)$.  However since $\ice$ is
necessarily unobservable, this formal identity has no operational
content.

\subsection{Effect of treatment on the treated}
\label{sec:ett}

Suppose that I am thinking of taking aspirin, and regard myself as
exchangeable with those individuals in the data who did in fact
receive aspirin --- though not necessarily with those who did not.  I
can then use the treated group to assess my hypothetical expected
response $\E(Y \cd F_T = 1)$ for $Y$, were I to take the aspirin; but
it seems I am not in a position to assess the contrasting hypothetical
expectation, $\E(Y \cd F_T = 0)$, and so cannot assess my personal
``effect of treatment''.  However, in the presence of a sufficient
covariate $U$---even if not observed---I may be able to do so.

We define the {\em effect of treated on the treated\/} as 
\begin{equation} 
  \label{eq:ettu} 
  \ett := \E(\sce_U \cd T =1, F_T = \idle). 
\end{equation}  
That is, $\ett$ is the average, in the observational regime, of the
specific causal effect (defined relative to $U$), over those
individuals who did in fact receive the asprin, $T=1$ --- and are thus
``like me''.

It might appear that, in the presence of a choice over which
sufficient covariate $U$ to use in \eqref{ettu}, that choice might
affect the value of $\ett$.  Fortunately it turns out that this is not
so, on account of the following result \cite{sgg/apd:ett}:
\begin{thm} 
  \label{thm:well}
  Suppose $\Pr(T=1 \cd F_T=\idle)>0$.  Then, for any sufficient
  covariate $U$, $\ett$ defined by \eqref{ettu} satisfies
  \begin{equation}
    \label{eq:ettwell} 
    \ett = \frac{\E(Y \cd F_T = \idle) - \E(Y \cd F_T = 0)}{\Pr(T=1 \cd F_T = \idle)}.     
  \end{equation} 
\end{thm} 
 
We have previously noted that, within the PR framework, we can
formally regard the pair $\bY$ of potential responses as a sufficient
covariate.  In that case the $\sce$ becomes the \ice, $Y_1 - Y_0$, and
\eqref{ettu} delivers $\ett = \E(Y_1 - Y_0 \cd T =1, F_T = \idle)$,
which is the usual PR definition of $\ett$.  However the above
argument shows that the PR framework is inessential for defining this
quantity.

Formula \eqref{ettwell} shows that we can identify $\ett$ whenever we
can observe the response $Y$ in the observational regime ($F_T =
\idle$), and also in a sample of people from whom the treatment was
withheld ($F_T = 0)$.  And although the definition of \ett\ supposes
the existence of some sufficient covariate, it is not necessary to
have observations on it.

\subsection{Reduction of sufficient covariate}
\label{sec:reduction}

Suppose $U$ is a sufficient covariate.  A function $V$ of $U$ is a
{\em sufficient reduction\/} of $U$ if $V$ is itself a sufficient
covariate.  Since property~\eqref{prop1} for $V$ follows immediately
from the same property for $U$, we only need investigate whether
property~\eqref{prop2v} holds for $V$:
\begin{equation}
  \label{eq:prop2v}
  \ind Y  {F_T} {(V,T)}.
\end{equation}

There are various additional conditions we can impose to ensure this.
One is the following:
\begin{cond}[Treatment-sufficient reduction]
  \label{cond:propred}
     \begin{equation}
      \label{eq:treatsuff0}
      \ind T U {(V, F_T=\idle)}.
    \end{equation}
    That is, in the observational regime, the choice of treatment
    depends on $U$ only through the value of $V$.
\end{cond}
Note that this condition does not involve the outcome variable $Y$ ---
except for the essential requirement that the starting variable $U$
itself be a sufficient covariate for the effect of $T$ on $Y$.  Also
note that, since $T$ is constant in any interventional regime,
\eqref{treatsuff0} is equivalent to
\begin{equation}
  \label{eq:treatsuff}
  \ind T U {(V, F_T)}.
\end{equation}
Also, since $V$ is a function of $U$, we trivially have
\begin{equation}
  \label{eq:star}
  \ind V {F_T} U,
\end{equation}
as well as 
\begin{equation}
  \label{eq:starstar}
  \ind Y V {(U, T, F_T)}.
\end{equation}

The following result now follows on applying the moralisation
criterion to the ID of \figref{treatsuff}, which faithfully represents
the conditional independence properties \eqref{prop1}, \eqref{star},
\eqref{treatsuff}, \eqref{prop2} and \eqref{starstar}, to deduce
\eqref{prop2v}:
\begin{figure}[htbp]
  \begin{center}
    \resizebox{2.5in}{!}{\includegraphics{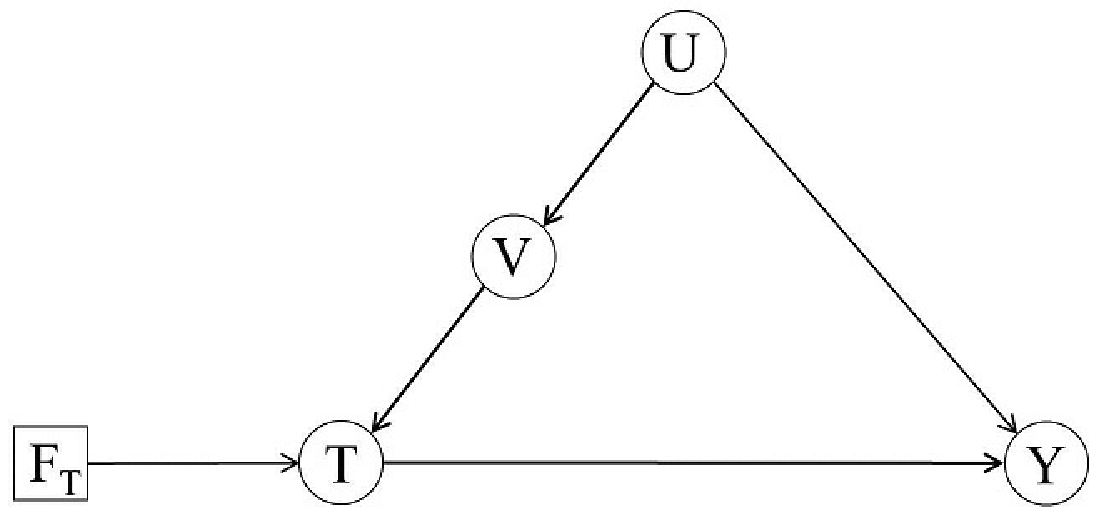}}
    \caption{Treatment-sufficient reduction}
    \label{fig:treatsuff}
  \end{center}
\end{figure}
\begin{thm}
  \label{thm:either}
  Suppose $U$ is a sufficient covariate, and let be $V$ be a function
  of $U$ such that \condref{propred} holds.  Then $V$ is a sufficient
  covariate.
\end{thm}
 
\subsubsection{Propensity score}
\label{sec:propensity}
An alternative description of treatment-sufficient reduction is as
follows.  Using P1, the defining property \eqref{treatsuff0} can be
expressed as
\begin{equation}
  \label{eq:balance}
  \ind U T{(V,F_T=\idle)}.
\end{equation}
In this form it asserts that, in the observational regime, the
conditional distribution of $U$ given $V$ is the same, whether further
conditioned on $T=0$, or on $T=1$: that is to say, $V$ is a {\em
  balancing score\/} for $U$ \cite{rosenbaum/rubin:propensity}.
Property \eqref{balance} can also be fruitfully interpreted as
follows.  Consider the family ${\cal Q} = \{Q_0, Q_1\}$ comprising the
pair of observational conditional distibutions for $U$, given,
respectively, $T=0$ and $T=1$.  Then \eqref{balance} asserts that $V$
is a {\em sufficient statistic\/} (in the usual Fisherian sense) for
this family.  In particular, a {\em minimal\/} treatment-sufficient
reduction is obtained as a minimal sufficient statistic for ${\cal
  Q}$: \viz, any $(1,1)$-function of the likelihood ratio statistic
$\Lambda := q_1(X)/q_0(X)$.  We might term such a minimal
treatment-sufficient covariate a {\em propensity variable\/}, since
one form for it is the treatment-assignment probability
\begin{equation}
  \label{eq:bayes}
  \Pi:  =\Pr(T=1 \cd U, F_T = \idle) = \pi\,\Lambda/(1-\pi + \pi\Lambda)
\end{equation}
(where $\pi := \Pr(T=1 \cd F_T = \idle)$), which is known as the {\em
  propensity score\/} \cite{rosenbaum/rubin:propensity}.  Either
$\Lambda$ or $\Pi$ supplies a $1$-dimensional sufficient reduction of
the orginal, perhaps highly multivariate, sufficient covariate
$U$.\footnote{However, this property may not be as useful as may first
  appear \cite{hg/apd:aistats2010}.}

\subsection{{\em do\/}-calculus}
\label{sec:docalc}

We here make use of the notation of \textcite{pearl:book} in which \eg
$p(y \cd x, \check z)$ refers to $\Pr(Y = y \cd X=x, F_Z = z)$, it
being implicit that $z \neq \idle$, and all unmentioned intervention
variables are idle.

Let $X$, $Y$, $Z$, $W$ be arbitrary sets of variables in a problem
also involving intervention variables.  The following rules follow
immediately from the definition of conditional
independence.\footnote{We assume throughout any positivity conditions
  required to ensure that the relevant conditional probabilities are
  well-defined.}

\begin{description}
\item[Rule 1 (Insertion/deletion of observations)] If $\ind Y Z {(X,
        F_X\neq \idle, W)}$ then
      \begin{equation}
        \label{eq:dorule1}
        p(y \cd \check x, z, w)= p(y \cd \check x, w).
      \end{equation}
\item[Rule 2 (Action/observation exchange)] If $\ind Y {F_Z} {(X,
        F_X\neq \idle, Z, W)}$, then
      \begin{equation}
        \label{eq:dorule2}
        p(y \cd \check x, \check z, w)= p(y \cd \check x, z, w).
      \end{equation}
\item[Rule 3 (Insertion/deletion of actions)] If $\ind Y {F_Z} {(X,
        F_X\neq \idle, W)}$, then
      \begin{equation}
        \label{eq:dorule3}
        p(y \cd \check x, \check z, w)= p(y \cd \check x, w).
      \end{equation}
\end{description}

Successive application of these rules, coupled with the property $F_X
= x \Rightarrow X=x$ and the laws of probability, can sometimes allow
one to express a ``causal'' expression in purely observational terms.
This was the essence of the argument in \secref{unconf} above, which
(assuming for simplicity that all variables are discrete) can be
expressed in general terms as:
\begin{theorem}[Back-door formula]
\label{thm:backdoor}
Suppose that
\begin{eqnarray}
  \label{eq:backdoorci1}
  Z &\cip& F_X\\
  \label{eq:backdoorci2}
  Y & \cip & F_X \,\,|\,\, (X, Z).
\end{eqnarray}
Then
\begin{equation}
  \label{eq:backdoorciform}
  p(y \cd \check x) =  \sum_z p(y\cd Z=z, X=x)\,p(Z=z).
\end{equation}
\end{theorem}

The most usual application of this {\em do-calculus\/} is for a model
represented by a Pearlian DAG.  However it is easiest to work with the
augmented DAG.\footnote{Pearl's analysis, like its precursor in
  \textcite{spirtes:ccd}, works with equivalent, somewhat more
  complex, formulations in terms of unaugmented DAGs.}  We first note
that conditioning on $F_X \neq \idle$ has the effect of removing all
arrows incoming to the set $X$ other than from $F_X$.  The resulting
reduced DAG can then be interrogated, using the usual moralisation
criterion, to deduce conditional independence properties that can be
used as input to Rules~1--3.\footnote{In fact Rule~1 is now
  redundant.}\@ In this context it can be shown constructively
\cite{shpitser:uai06,shpitser:ncai06,huang:uai06} that, whenever there
exists a reduction of a causal expression to purely observational
terms, it can be found by applying the {\em do\/}-calculus.

\section{Instrumental Variables}
\label{sec:instrument}
\begin{figure}[htb]
  \begin{center}
    \resizebox{1.5in}{!}{\includegraphics{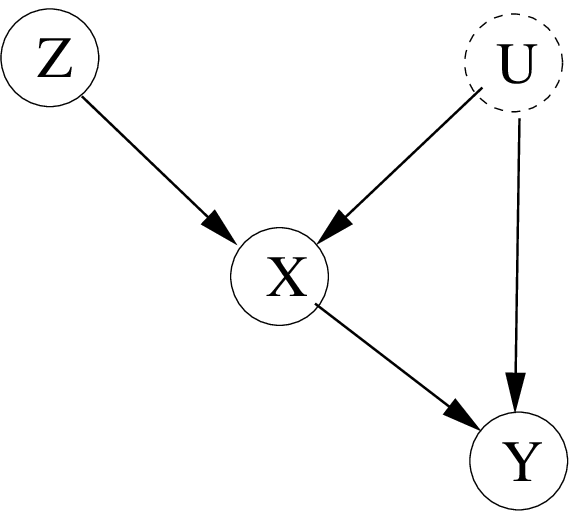}}
    \caption{Instrumental variable}
    \label{fig:instrument0}
  \end{center}
\end{figure}

In the presence of an unobserved sufficient covariate $U$, it is
typically not possible to estimate the average causal effect, \ace, of
a treatment variable $X$ on a response variable $Y$ from observational
data.  Some progress can be made if we can assume the existence of an
observable {\em instrumental variable\/} $Z$, which can be thought of
as an imperfect proxy for an intervention.  The assumptions required
in such a case are typically expressed in informal terms such as
\cite{martens:06}:
\begin{enumerate}
\item \label{it:martens1} $Z$ has a {causal} {effect} on $X$
\item \label{it:martens2} $Z$ {affects} the outcome $Y$ {only}
  {through} $X$ (``{no} {direct} {effect} of $Z$ on $Y$'')
\item \label{it:martens3} $Z$ does not share {common} {causes} with
  the outcome $Y$ (``{no} {confounding} of the {effect} of $Z$ on
  $Y$'').
\end{enumerate}
These might be formalised as observational conditional independence
properties, such as:
\begin{eqnarray}
  \label{eq:obsci0}
  X &\not\!\!\!\cip& Z\\
  \label{eq:obsci2}
  U &\cip&  Z\\
  \label{eq:obsci4}
  Y &\cip& Z \,\cd\, (X, U).
\end{eqnarray}
Note the analogy between \eqref{obsci2} and \eqref{prop1}, and
\eqref{obsci4} and \eqref{prop2}, where (with $T$ relabelled as $X$)
$Z$ takes the place of $F_X$.  However, unlike the case of an imposed
intervention, $Z$ does not determine the value of $X$, but merely has
some association with it, as described by \eqref{obsci0}.  These
assumptions are represented by the DAG of
\figref{instrument0}.\footnote{For \eqref{obsci0}, we need to assume
  that this is a faithful representation.}

\begin{figure}[htb]
  \begin{center}
    \resizebox{1.5in}{!}{\includegraphics{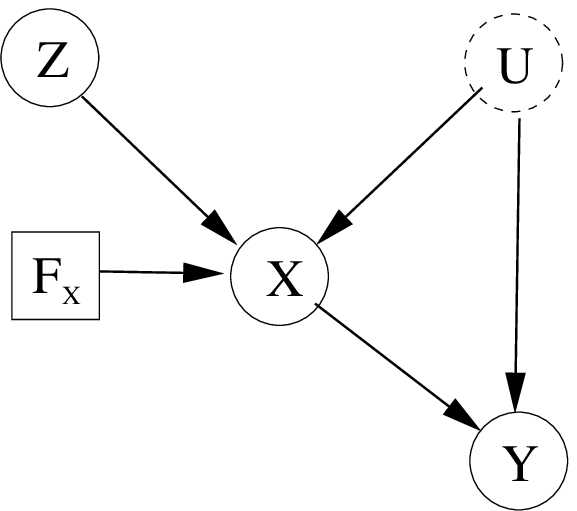}}
    \caption{Instrumental variable with regimes}
    \label{fig:instrument}
  \end{center}
\end{figure}

For all that this might be a fruitful analogy, requirements
\eqref{obsci0}--\eqref{obsci4}, and \figref{instrument0}, leave
something to be desired: since they relate solely to the observational
regime, they can not, of themselves, have any causal consequences ---
at best these are left implicit, which leaves room for confusion.  It
is far better to make the requisite {\em causal\/} assumptions
explicit.  We do this by elaborating \figref{instrument0} to
explicitly include the nonstochastic regime indicator $F_X$ for $X$,
as in \figref{instrument}.
For $F_X = \idle$ this recovers the assumptions encoded in
\figref{instrument0}; but in addition it relates the observational
structure to what would happen under an intervention to set $X$.  In
particular, it clarifies that $U$ is assumed to be a sufficient
covariate for the effect of $X$ on $Y$, and further
encodes:\footnote{\figref{instrument} also encodes the additional, but
  inessential, property $\indo Z {F_X}$.}
\begin{eqnarray}
  \label{eq:instr1}
  U &\cip& Z  \,\,\cd \,\, F_X\\
  \label{eq:instr2}
   Y &\cip& Z \,\,\cd \,\,(X, U, F_X).
\end{eqnarray}
Properties~\eqref{instr2} and \eqref{instr1} extend \eqref{obsci2} and
\eqref{obsci4} to apply under intervention, as well as
observationally.

\subsection{Linear model}
\label{sec:linmod}

Suppose now all the observables are univariate, and we can describe
the dependence of $Y$ on $(X, U)$ (which we have assumed the same in
all regimes) by a linear model:
\begin{equation}
  \label{eq:linmod}
  \E(Y \cd X, U, F_X) = W + \beta\,X
\end{equation}
for some function $W$ of $U$.  

We deduce
\begin{displaymath}
  \E( Y \cd F_X=x) = w_0 + \beta\, x,
\end{displaymath}
where $w_0 := \E(W \cd F_X = x)$ is a constant independent of $x$,
since $\indo U {F_X}$.  Thus $\beta$ can be interpreted causally, as
describing how the mean of $Y$ changes in response to manipulation of
$X$.  Our aim is to identify $\beta$.

By \eqref{instr2}, \eqref{linmod} is also $\E(Y \cd X, Z, U, F_X =
\idle)$.  Then
\begin{displaymath}
  \E(Y \cd Z, F_X = \idle) = \E(W \cd Z, F_X =\idle) + \beta\,\E(X
  \cd Z, F_X =\idle).
\end{displaymath}  
But by \eqref{instr1} the first term on the right-hand side is
constant.  Thus
\begin{equation}
  \label{eq:instrform}
  \E(Y \cd Z, F_X =\idle) =  \mbox{constant} + \beta\,\E(X \cd Z,
  F_X =\idle).
\end{equation}
Equation~\eqref{instrform} relates two functions of $Z$, each of which
can be identified from observational data.\label{page:fninstr}
Consequently (so long as neither side is constant) we can identify the
causal parameter $\beta$ from such data.  Indeed it readily follows
from \eqref{instrform} that (in the observational regime) $\beta =
\cov(Y,Z)/\cov(X,Z)$, which can be estimated by the ratio of the
coefficients of $Z$ in the sample linear regressions of $Y$ on $Z$ and
of $X$ on $Z$.
 
\subsection{Binary variables}
\label{sec:binary}

When all the observable variables $Z,X, Y$ are binary, without making
further assumptions we can not fully identify the ``causal
probability'' $P(Y=1 \cd F_X =x)$ from observational data.  However,
we can develop inequalities it must satisfy.  This approach was
instigated by \textcite{manski}.  His inequalities were refined by
\textcite{b+p:jasa}, under the strong additional condition of
deterministic dependence\footnote{An alternative interpretation of
  this condition is in terms of potential outcomes.}  of $X$ on
$(Z,U)$ and of $Y$ on $(X,U)$.  This condition was shown to be
unnecessary by \textcite{apd:hsss}, where a fully stochastic
decision-theoretic approach was developed.  In either approach, the
analysis involves subtle convex duality arguments.

\section{Dynamic Treatment Strategies}
\label{sec:dynamic}
\begin{figure}[h]
  \begin{center}
    \resizebox{2.8in}{!}{\includegraphics{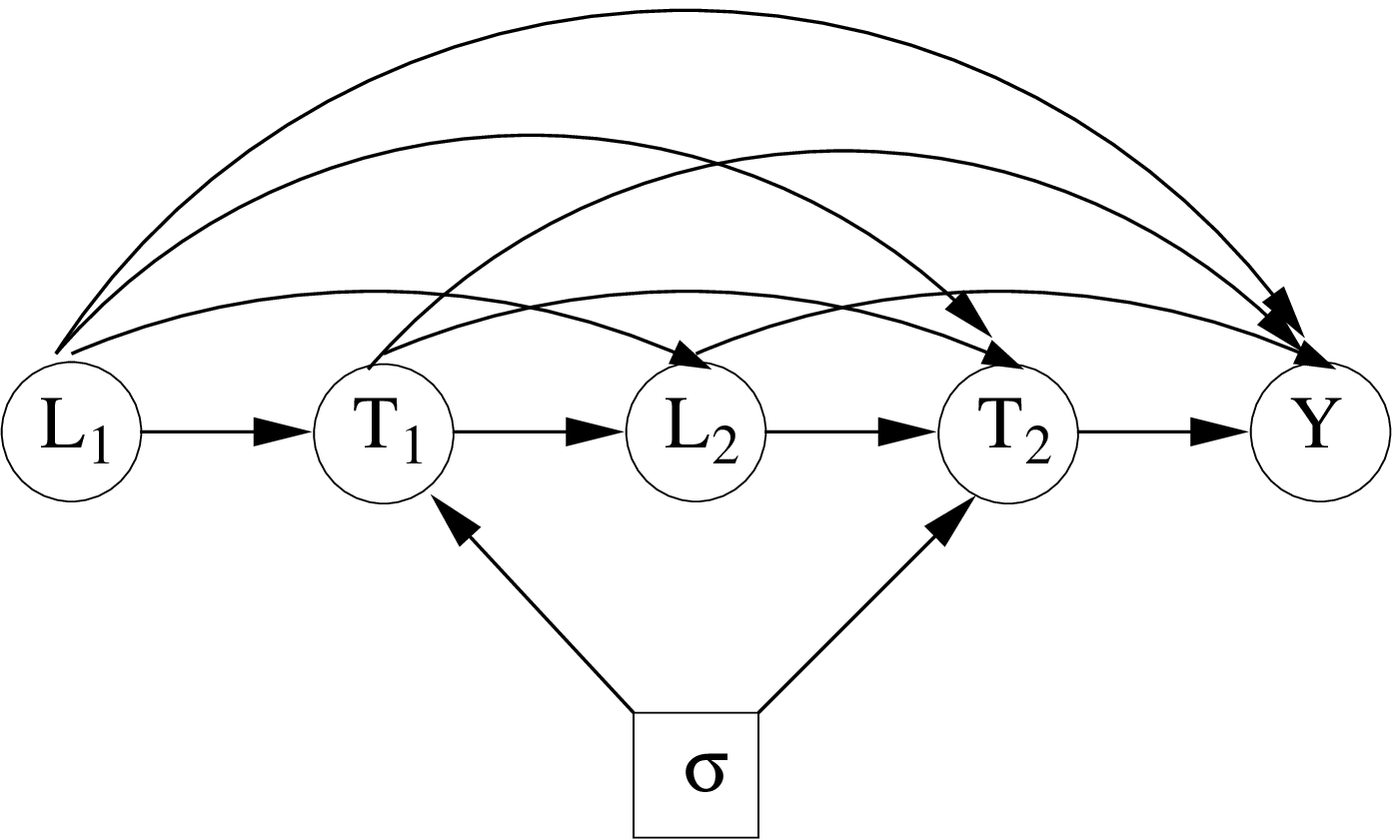}}
    \caption{Sequential ignorability}
    \label{fig:seqignor}
  \end{center}
\end{figure}
In the ID of \figref{seqignor}, the $L$'s represent attributes of a
patient, the $T$'s treatments that can be applied, and $Y$ a response
of interest.  These variables are supposed generated in the order
shown, each in response to all its predecessors.  The non-stochastic
regime indicator node $\sigma$ can take value \idle, indicating the
observational regime; otherwise, a value $\sigma=s$ describes a
hypothetical {\em treatment strategy\/}, specifying how treatment
$T_1$ should be chosen in response to observation of $L_1$, and how
$T_2$ should be chosen in response to observation of $(L_1, T_1,
L_2)$.  Typically such a strategy will prescribe deterministic
choices, but there is no difficulty in allowing further randomisation.
The task is to infer the {\em consequence\/}, $\E(Y \cd \sigma=s)$, of
such a hypothetical strategy from properties of the observational
regime $\sigma=\idle$.

\figref{seqignor} encodes the following conditional independencies:
\begin{eqnarray}
  \label{eq:L1}
  L_1 &\cip& \sigma\\
  \label{eq:L2}
  L_2 &\cip& \sigma \cd (L_1, T_1)\\
  \label{eq:Y}
  Y &\cip& \sigma \cd (L_1, T_1, L_2, T_2).
\end{eqnarray}
Condition~\eqref{Y}, for example, says that the distribution of $Y$,
given the previous variables $(L_1, T_1, L_2,T_2)$, in the
observational regime $\sigma=\idle$ would also apply under the
operation of an imposed strategy $\sigma=s$.  This is a ``no residual
confounding'' type of assumption, that might or might not be
appropriate.  When \eqref{L1}--\eqref{Y} apply, we say we have {\em
  sequential ignorability\/}.

We will always have
\begin{eqnarray}
  \label{eq:1}
  p(l_1, t_1, l_2, t_2, y\cd \sigma = s)& = & p(l_1 \cd \sigma = s)\\
  \label{eq:2}
  && { }\times p(t_1 \cd l_1, \sigma = s)\\
  \label{eq:3}
  && { }\times  p(l_2 \cd l_1, t_1, \sigma = s)\\
  \label{eq:4}
  && { }\times  p(t_2 \cd l_1, t_1, l_2, \sigma = s)\\
  \label{eq:5}
  && { }\times p(y \cd l_1, t_1, l_2, t_2, \sigma = s).
\end{eqnarray}
Now \eqref{2} and \eqref{4} are specified by the strategy $s$.  Also,
under sequential ignorability, in \eqref{1}, \eqref{3} and \eqref{5}
we can replace $\sigma=s$ by $\sigma = \idle$, so that those terms are
estimable from observational data.  We will thus have all the
ingredients needed to identify the joint distribution of all variables
under the strategy $\sigma=s$, and then by marginalisation we can
identify the desired consequence, $\E(Y \cd \sigma=s)$.  This
computation, which can be effectively restructured as a recursion
\cite{dd:ss}, reduces to the {\em $g$-computation\/} formula of
\textcite{jr:mm}.  That paper (see also \textcite{murphy:annrev}) set
the problem up in a PR framework, assuming the simultaneous existence
of potential responses $(L_{1s}, L_{2s}, Y_s)$ for each possible
strategy $\sigma = s$, subject to certain consistency requirements,
sequential ignorability then being expressed as a conditional
independence property involving these potential responses.  Our DT
approach is more straightforward to interpret, justify and implement,
as well as allowing for randomised strategies.

It will often be unrealistic to impose the ``no residual confounding''
assumptions of sequential ignorability, at least without further
justification.  Such an assumption might become more reasonable when
additional variables are added to the system: variables that could
not, however, be usable by the considered strategy $\sigma=s$.  In
such a case it is possible to add further conditions, generalising
those of \secref{nonconf}, which when acceptable would imply that we
will indeed have sequential ignorability.  For further details see
\textcite{dd:ss,apd/pc:2014}.

\section{Discussion}
\label{sec:disc}
The decision-theoretic language for causality has sometimes been
criticised for not being as rich as that of alternative approaches,
such as PR models, which can make statements, in their own
mathematical terms, that simply have no DT counterpart.  I regard this
as a strength, not a weakness: formal mathematical expressions (for
example, the variance of the \ice---see \secref{prmod}) that do not
relate directly to features of the real world are at best unnecessary,
and at worst dangerously misleading.  Within DT we are not plagued
with ``the fundamental problem of causal inference'' \cite{pwh:jasa},
which is only a self-created problem of the PR approach.  The DT
approach also fosters healthy scepticism of other methods, such as
``principal stratification'' \cite{frang/dbr:princstrat}, that depend
crucially on the philosophically perplexing assumption of the real
simultaneous existence of potential response pairs
\cite{apd/vd:princstrat}, together with necessarily untestable
assumptions about their properties.  Within the ambit of problems that
are well-posed, the DT framework has all the expressive power
necessary, uncluttered by unnecessary and distracting formal
mathematical ingredients.

\nocite{glymour/cooper:99}

\bibliographystyle{oupvar} \bibliography{strings,causal,allclean,ci}
\end{document}